\documentstyle{amsppt}

\magnification=1100 \NoBlackBoxes \nologo \hsize=14.5cm
\vsize=19.5cm

\def\span{\text{span}}
\def\even{\text{even}}
\def\odd{\text{odd}}

\def\inv{^{-1}}
\def\?{{\bf{??}}}

\def\Ndr{N_d^{\text{red}}}

\def\A{\Bbb A}

\def\F{\Cal F}

\def\C{\Bbb C}
\def\N{\Bbb N}

\def\P{\Bbb P}

\def\ls{\vskip.25in}
\def\ss{\vskip.15in}

\def\[{\big[}
\def\]{\big]}
\def\V{\bar{V}}

\def\O{\Cal O}

\def\Sym{\text{Sym}}

\def\rk{\text{rk}}

\def\m{\frak m}
\def\k{\underline{ k}}

\def\1/2{\frac{1}{2}}

\def\D{\frak D}

\def\I{\Cal I}

\def\im{\text{im}}

\def\2{{[2]}}
\def\l{\ell}
\topmatter
\title Normal bundles of rational curves in projective spaces
\endtitle
\author
Ziv Ran
\endauthor
\thanks{\raggedright{
Partially supported by NSA Grant
MDA904-02-1-0094}\newline\indent\raggedright{\ \ Updated version
at $\underline{www.math.ucr.edu/\tilde{\ } ziv/papers/}$} }
\endthanks
\date September 30, 2003\enddate

\address University of California, Riverside\endaddress
\email ziv\@math.ucr.edu\endemail
\rightheadtext { Normal Bundles} \leftheadtext{Z. Ran}

\abstract We determine the splitting
(isomorphism) type of the normal bundle of a
generic genus-0 curve with 1 or 2 components in
$\P^n,$ as well as the way the bundle deforms
locally with a general deformation of the curve.
We deduce an enumerative formula for divisorial
loci of smooth rational curves whose normal
bundle is of non-generic splitting
type.\endabstract
\endtopmatter
\comment ignored text\endcomment

\document
Rational curves in projective space, being essentially the same
thing as finite-dimensional vector spaces of rational functions in
one variable, are among the most elementary and classical objects
in Algebraic Geometry. In recent years it has become clear that
suitable (compact) parameter spaces, say $R_{n,d}$, for rational
curves of given degree $d$ in $\P^n$, are of fundamental
importance. Now the geometry of a moduli or parameter space like
$R_{n,d}$ is closely related to 'modular' subvarieties, i.e. ones
defined in terms of the family of curves (or other objects) that
it parametrizes. Now there are various ways of defining modular
subvarieties of $R_{n,d}$, for instance {\it{incidence}}
subvarieties, parametrizing curves incident to a given cycle in
$\P^n.$ Another such way involves vector bundles. Namely, given a
'reasonable function' $\Phi$ assigning to a curve $C\in R_{n,d}$ a
vector bundle $E_C$ on $C$, a theorem of Grothendieck says we have
a decomposition
$$E_C\simeq\bigoplus\O_C(k_i), k_1\geq k_2,...$$
where $\O_C(k)$ denotes the unique line bundle of degree $k$ on
$C$. The sequence $k.=k.(C)$, which is uniquely determined and
called the {\it{splitting type}} of $E_C$, varies
upper-semicontinuously, in an obvious sense, in terms of the
vector bundle and hence for a good function $\Phi$ we get a
stratification $R_{n,d}^\Phi(k.)$ of $R_{n,d}$ where the strata
consist of the curves $C$ with given sequence $k.(C)$.\par One way
to define an interesting, and reasonable, function $\Phi$ is to
fix a vector bundle $E$ on $\P^n$ and to set
$$E_C=E|_C.$$
The resulting stratification was studied in \cite{R5} where we
computed enumeratively its divisorial stratum. The main result of
this paper is an analogous computation in the case where $\Phi$ is
the 'normal bundle function', which assigns to a curve $C$ its
{\it{normal bundle}}
$$N_C=N_{C/\P^n}.$$ The splitting type $k.(C)$ of $N_C$, which we
call the {\it{normal type}} of $C$, is a natural global numerical
invariant of the embedding $C\subset\P^n$, perhaps the most
fundamental such invariant beyond the degree, and thus the problem
of enumerating curves with given normal type seems a natural one.
Despite the existence of a natural map from the restricted tangent
bundle $$T_{\P^n}|_C\to N_C$$ there is, in fact, little
relationship between the splitting types of these bundles and
consequently, though the main enumerative result of this paper
(Theorem 9.1) is an analogue for normal bundles of the result of
\cite{R5} for restricted bundles, there is in reality little of
substance in common between the two papers.\par In very broad
outline, the proof of Theorem 9.1 goes as follows. To begin with,
for the enumerative question to make sense it is necessary that
the normal type $k.(C)$ of the generic curve in $R_{d,n}$ be 'well
behaved'. This turns out to mean that this type is {\it{almost
balanced}} in the sense that
$$k_1(C)-k_{n-1}(C)\leq 1.$$ Assuming this, we need that the
'jump' of $k.(C)$ (from its generic value) occurs  in codimension
1. For that, it turns out that a necessary and sufficient
condition is that $(d,n)$ be a {\it{perfect pair}} in the sense
that
$$(n-1)|2d-2,$$ or equivalently, that the type $k.(C)$ of the
generic $C$ should be {\it{balanced}}, in the sense that
$k_1(C)=k_{n-1}(C)$. Fixing a perfect pair $(d,n)$, we may, as in
\cite{R5}, consider a generic {\it{incidence pencil}}, i.e. an
($\infty^1$) $B$ of curves in $R_{n,d}$ defined by incidence
conditions, and the (smooth) surface $X$ swept out by the curves
in $B$, and on $X$ an appropriate vector bundle $G$ that is a
twist of the 'relative lci normal bundle'. Then for smooth members
of the pencil, unbalanced type can be interpreted in terms of the
cohomology of $G$, as the local length of a suitable $R^1$ sheaf,
which is of finite support and length, and where the other $R^i$
vanish. Now the Grothendieck-Riemann-Roch formula gives an
expression for the total length of $R^1$. Therefore to complete
the proof it 'remains' to evaluate the contribution from singular,
i.e. reducible fibres.\par This evaluation turns out, in reality,
to be possibly the most involved part of the story. For, in
contrast with the case of restricted bundles, the (lci) normal
bundle of a reducible curve is quite often ill behaved; viz. there
are natural notions of balanced and almost balanced for bundles on
(2-component) reducible curves, and the normal bundle to a curve
with a degenerate component is usually {\it{not}} almost balanced
(and even in cases when the latter bundle is almost balanced, that
fact is relatively subtle to prove). The unbalancedness has some
significant implications. First, it makes substantially more
difficult an attempt to prove generic almost balancedness by
specialization; second, and more consequentially, it means that in
a pencil as above it is not sufficient to compute the cohomology
{\it{on}} reducible fibres, but one must compute it in a
neighborhood as well, i.e. compute the (length of) the entire
local $R^1$ module, a substantially more involved computation
which ultimately shows that $R^1$ is 'as small as possible', given
the $H^1$ on the reducible fibre, i.e. that $R^1$ is killed by the
maximal ideal. This result, that we call 'cohomological
quasitransversality', is established by constructing an explicit
smoothing of 'binomial' type for which the requisite property
reduces to a combinatorial property of the exponents that we
establish by a somewhat drawn-out, but quite elementary
argument.\par The paper is organized as follows. In \S1 we discuss
elementary modifications of vector bundles and some of their basic
properties. In \S2 we discuss the (lci) relative normal bundle in
a family of curves in a smooth ambient variety. In \S3 we give a
preliminary elementary discussion of normal bundles to rational
curves, especially the rational normal curve, in projective space.
In \S4 we study in further detail the normal bundle to a rational
normal curve and give a geometric interpretation of its splitting.
In \S5 we give a general elementary discussion of vector bundles
on rational trees, especially 2-component trees that we call
{\it{rational angles}}. On a rational angle, every vector bundle
is a direct sum of lines bundles, but this fails on a general
rational tree by Example 5.6. We study especially almost balanced
bundles and their deformations and specializations. In \S6 we give
a complete determination of the normal bundles of general rational
curves and rational angles in $\P^n.$ We find that the normal
bundle is almost balanced for general rational curves of degree
$d\geq n$; for a general rational angle $C_a\cup C_b $ we find
that the normal bundle is almost balanced if both $a,b\geq n$ but
usually not otherwise. \S7 is preparatory to \S8 in which we prove
the cohomological quasitransversality result mentioned above.
After all these preparations, the proof of the main enumerative
result, Theorem 9.1, is a straightforward adaptation of that of
the main result of \cite{R5}, and like it uses the intersection
calculus on incidence pencils, developed in earlier papers and
reviewed in an Appendix.\par An interesting question not settled
by our work is that of irreducibility of the locus of curves of
degree $d$ in $\P^n, n\geq 4,$ with given normal type $(k.)$. This
seems open even in case $(d,n)$ is a perfect pair, so that the
locus in question is of codimension 1. For $n=3$ irreducibility is
known by \cite{EV2}. By contrast, the analogous irreducibility
property for the restricted tangent bundle holds trivially,
because, by the Euler sequence, the locus in question is
parametrized by an open subset of the vector space
$$\text{Ext}^1(\bigoplus\O_{\P^1}(k_i),
\O_{\P^1})=H^1(\bigoplus\O_{\P^1}(-k_i))$$ (see
also \cite{Ram}).

\subheading{Acknowledgment} Some of this work was carried out
while I was visiting the IHES, and I would like to thank the
Institute and its staff for providing an ideal working
environment. I also thank the organizers of the conference
'Current Geometry', held in Naples in June 2003, for an
opportunity to present this work. \subheading{Notational
conventions} On a nonsingular rational curve $C$ we will denote by
$\O_C(k)$ or $\O(k)$ the unique line bundle of degree $k$ on $C$.
On a reducible curve, we will use notation like $\O(j\cup k)$ to
denote the like bundle having degrees $j,k$ on the respective
components. For any sheaf or abelian group $L$ and natural number
$m$, $mL$ usually denotes $\bigoplus\limits_1^m L$ (unless it's
clear that $L$ is being viewed as a divisor, in which case $mL$ is
its multiple as such). We will be working over the groundfield
$k=\C$, and using only closed points $p$, so that the residue
field $k(p)\simeq\C$ always, nonetheless the notation $k(p)$ or
$\underline{k}(p)$ will be used. For a coherent sheaf $L$ and a
point $p$, the fibre $L(p):=L\otimes k(p).$
 \subheading{1. Elementary Modifications}\ss Let $X$ be a
reduced algebraic scheme and $E$ a locally free
coherent sheaf on $X$. Consider a quotient of the
form
$$E\overset\phi\to\to q\to 0 \tag 1.1$$
where $q$ is a locally free $\O_\Sigma$-module for a Cartier
divisor $\Sigma$ on $X$. It is easy to see that the kernel
$E'\subset E$ of $\phi$ is a locally free sheaf called the
{\it{elementary reduction}} of $E$ corresponding to $q$, denoted
$M(E,q)$ (of course in reality, $E'$ depends on $\phi$ and not
just $q$). We will primarily be interested in the cases \par

(1) X is a curve, $\Sigma$ is a reduced smooth point on $X$ and
$q$ has length 1;\par

(2) X is a smooth surface and $\Sigma$ is a reduced curve on $X$.
\par\noindent
From the defining exact sequence
$$0\to M(E,q)\to E\to q\to 0\tag 1.2$$
it is easy to compute the Chern classes of $M(E,q)$; e.g. in case
(1) we get
$$c_1(M(E,q))=c_1(E)-\l(q).$$
Dualizing (1.2), we obtain the exact sequence
$$0\to E^*\to M(E,q)^*\to \sigma^*\to 0$$
where $\sigma$ (and $S$ below) are defined by the exact diagram
$$\matrix 0\to \sigma\to &M(E,q)\otimes\O_\Sigma&
\to &E\otimes\O_\Sigma&\to q\to 0\\
&\searrow&&\nearrow&\\
&&S&&\\
&\nearrow&&\searrow&\\
0&&&&0\endmatrix\tag 1.3$$ and $\sigma^*$ means $\O_\Sigma-$dual.
Thus
$$M(M(E,q)^*, \sigma^*)^*=E.\tag 1.4$$
The bundle $M(E,q)^*$ may be called the {\it{elementary
enlargement}} of $E^*$ corresponding to the subbundle
$q^*\subseteq E^*\otimes\O_\Sigma.$ An {\it{elementary
modification}} is an elementary reduction or enlargement. It is
also easy to see from (1.3) that
$$M(M(E,q),S)=E(-\Sigma).\tag 1.5$$
Consequently, an elementary reduction of $E$ is an elementary
enlargement of $E(-\Sigma)$. In practice this means that it
suffices to work with elementary reductions, which are more
convenient than enlargements.
\par
Now suppose we have an exact sequence of vector bundles
$$ 0\to F\to E\to G\to 0.$$
We will say that $F$ {\it{survives}} (resp. {\it{gets chopped}})
in the elementary modification (1.1) if the induced map
$$F\to q$$
is zero (resp. surjective). If $F$ survives, then
considering $F$ as subsheaf of $E$, we have
$$F\otimes\O_\Sigma\subset S$$
while, considering $F$ as subsheaf of $E'$, we have
$$F\otimes\O_\Sigma\cap \sigma =0,$$
and we have an exact sequence
$$0\to F\to E'\to G'\to 0$$
where $G'=M (G,q).$ If $F$ gets chopped, we get a subbundle (i.e.
locally split subsheaf)
$$F':=M(F,q)\to E'$$
and
$$\sigma\subseteq F'\otimes\O_\Sigma ,$$
and we have an exact sequence
$$0\to F'\to E'\to G\to 0.$$

\subheading{2. The lci normal bundle} \ss Let $\pi:X\to B$ be a
family of nodal curves and $f:X\to Y$ a generically 1-1 map to a
smooth variety. Assume that $f$ is unramified on all fibres and an
embedding on almost all fibres, including all singular ones. We
have an exact sequence
$$f^*\Omega_Y\overset{df}\to{\to}\omega_{X/B}\to
q_\Sigma\to 0\tag 2.1$$ where $\Sigma$ is the critical locus, i.e.
the locus of singular points of fibres of $\pi$, and $q_\Sigma$ is
a skyscraper sheaf with length 1 at each point of $\Sigma$. We
denote $\ker(df)$  by $N^*_{f/B}$ or just $N^*$ for short, so that
we have a basic exact sequence
$$0\to N^*\to f^*\Omega_Y\overset{df}\to{\to}\omega_{X/B}\to
q_\Sigma\to 0,\tag 2.2$$ where the image of $df$ coincides with
$\I_\Sigma\omega_{X/B}(L)$. In light of the exact sequence
$$0\to
N^*(L)\to f^*\Omega_Y(L)\overset{df}\to{\to}\Omega_{X/B}(L)\to
0,$$ we see that $\Omega_{X/B}=\I_\Sigma\omega_{X/B}$, and in
particular $\Omega_{X/B}$ is torsion-free and $\omega_{X/B}$ is
its double dual.

For any line bundle $L_Y$ on $Y$ and $L=f^*L_Y$, note the  exact
diagram%
$$\matrix &&&0&&0\\&&&\downarrow&&\downarrow\\
0\to
&N^*(L)&\to&f^*\Omega_Y(L)&\overset{df}\to{\to}&\Omega_{X/B}(L)&\to&
0 \\
&\parallel&&\downarrow&&\downarrow\\
0\to&N^*(L)&\to&f^*P^1_Y(L_Y)&\to&P^1_{X/B}(L)&\to 0\\
&&&\downarrow&&\downarrow\\
&&&L&=&L\\
&&&\downarrow&&\downarrow\\
&&&0&&0\\
\endmatrix \tag 2.3$$

Let $P^+_{X/B}(L)$ be the double dual of $P^1_{X/B}(L)$, which may
be called the sheaf of 'relative dualizing principal parts' of
$L$, and which is clearly locally free and fits in a diagram
$$\matrix 0\to&\Omega_{X/B}(L)&\to&P^1_{X/B}(L)&\to&L&\to 0\\
&\downarrow&&\downarrow&&\parallel\\
0\to&\omega_{X/B}(L)&\to&P^+_{X/B}(L)&\to&L&\to 0.\endmatrix\tag
2.4$$ Then we get an exact sequence
$$0\to N^*(L)\to f^*(P^1_Y(L_Y))\to P^+_{X/B}(L)\to q_\Sigma(L)\to
0\tag 2.5$$ This sequence is especially useful when $Y=\P^n=\P(V)$
and $L_Y=\O(1)$-the case of principal interest to us- in which we
have, as is well known $P^1_Y(L_Y)=V^*\otimes\O_Y.$

 Now by easy and well known local
computations (partly reproduced below), $N^*$ is a locally free
sheaf of rank $n-1:=\dim(Y)-1$ and is called the {\it{(relative)
lci conormal bundle}} of the map $f$. For any fibre $X_b$ we have,
setting $C=f(X_b)$,
$$N^*_{X_b}\simeq f^*(\I_{C}/\I^2_C)/{\text{(torsion)}}.$$
In particular, if $f|_{X_b}$ is an embedding, then

$$N^*_{X_b}\simeq f^*(\I_{C}/\I^2_C);$$
generally, if $x\in X_b$ is a smooth point then the fibre of $N^*$
at $x$, denoted $N^*(x)$, is canonically isomorphic to the
conormal space at $f(x)$ of the unique branch of $C$ coming from
an analytic neighborhood of $x$ on $X_b.$\par Let us analyze the
situation locally at a fibre node $p\in\Sigma.$ For simplicity we
assume the fibre $X_0$ through $p$ is a union of 2 smooth fibres
$X_1, X_2$- this is the case we will need. We may choose local
coordinates $x_1,...,x_n$ on $Y$ so that $X_i$ maps to the
$x_i$-axis, $i=1,2$ so $C=f(X_0)$ is locally defined by
$$x_1x_2=x_3=\cdots =x_n=0.$$
Let $$C_i=f(X_i), N^*_i=f^*(\I_{C_1}/I^2_{C_1}), i=1,2.$$
Then we have an exact sequence
$$0\to
N^*|_{X_1}\to f^*(\I_{C_1}/I^2_{C_1})\to q_p\to 0$$
which restricts at $p$ to
$$0\to
<\sigma_p>\to N^*(p)\to
N_1^*(p)\to q_p\to 0.$$
In terms of bases, the latter sequence can be written
$$0\to <x_1x_2>\to <x_1x_2,x_3,\cdots , x_n>\to$$$$
<x_2,\cdots , x_n>\to <x_2,\cdots , x_n>/<x_3,\cdots ,x_n>
\to 0.$$
We usually set $M=N^* (p)$ and $\sigma=\sigma_p$ is called
the {\it{singular}} element of $M$ (well-defined up to
scalar).
Note that the image $S$
of $N^*(\sigma)\to N_1^*(\sigma)$, i.e.
$<x_3,\cdots,x_n>$, is just the Zariski conormal space
to $C$ at $f(p)$, and we have an exact sequence
$$0\to <\sigma>\to M\to S\to 0.\tag 2.6$$
\remark{Example 2.1} Let $C=L_1\cup L_2$ where
$L_1,L_2\subset\P^n$ are distinct lines meeting
at $p.$ Then
$$N^*_C|_{L_1}=M(N^*_{L_1},q)$$
where $q$ is the quotient of $N^*_{L_1}$
corresponding to the tangent direction of $L_2$
at $p$. Since
$$N^*_{L_1}\simeq N^*_{L_1}(\sigma)\otimes\O_{L_1}(-1),$$
it has a unique surviving subsheaf isomorphic to
$(n-2)\O(-1)$, viz. $S\otimes\O(-1)$ (in the obvious
sense). Choosing any chopped subsheaf isomorphic to
$\O(-1),$
hence complementary to $S\otimes\O_{L_1}(-1)$,
 we get a subbundle
$$\O_{L_1}(-2)\subset N^*_C|_{L_1}$$
whose fibre at $p$ is $q.$ Since the fibre at $p$
 of the
above surviving subsheaf $S\otimes\O_{L_1}(-1)$,
as subsheaf of $N^*_C|_{L_1}$, is the Zariski
conormal space, it is complementary to $q$ at
$p$, so we get a splitting
$$N^*_C|_{L_1}\simeq (n-2)\O_{L_1}(-1)\oplus \O_{L_1}(-2)$$
where the $(n-2)\O_{L_1}(-1)$ summand, and the
fibre of the $\O_{L_1}(-2)$ summand at $p$ (but
not the summand itself), are uniquely determined,
the latter being $q$. As already observed,
 the fibre at $p$
of the $(n-2)\O_{L_1}(-1)$ summand is just the Zariski conormal
space which, in local coordinates as above, has $x_3,...,x_n$ as
basis.\par Since $N^*|_{L_2}$ can be analyzed similarly, we find
that the $(n-2)\O_{L_1}(-1)$ and $(n-2)\O_{L_2}(-1)$ glue at $p$
as do the $\O_{L_1}(-2)$ and $\O_{L_2}(-2)$ summands, so we get a
splitting
$$N^*_C=(n-2)\O((-1)\cup(-1))\oplus\O((-2)\cup(-2))$$
where $\O(a\cup b)$  is the line bundle on $C$ having degree $a$
on $L_1$ and $b$ on $L_2$ (such notation will be used throughout
the paper). This result if of course obvious from the fact that
$C$ is a $(2,1^{n-2})$ complete intersection, but is nonetheless
enlightening in that it shows that the two positive subhseaves
$$(n-2)\O_{L_i}(-1)\subset N^*_C|_{L_i}, i=1,2$$
are {\it{not}} mutually in general position at $p$, contrary to
what one might naively have expected. This contrasts with the
situation when normal bundles are replaced by restrictions of a
fixed (and suitable) bundle on $\P^n$, for instance the tangent
bundle (cf \cite{R5}, \S2).

\endremark
To formalize the sort of situation typically
encountered in analyzing the normal bundle of a
reducible curve, it is convenient to introduce
some definitions. \proclaim{Definition 2.2} Let
$M$ be a vector space with a distinguished
1-dimensional subspace $<\sigma>$. A pair $F.,
G.$   of increasing filtrations on $M$ are said
to be in {\it{relative general position}} (with
respect to $\sigma$, if that is not understood)
provided we have for each $i,j$ that \item{(i)}
whenever $\sigma\not\in F_i\cap G_j,$ $F_i$ and
$G_j$ are in general position in $M$; \item
{(ii)} whenever $\sigma\in F_i\cap G_j,$
$F_i/\sigma, G_j/\sigma$ are in general position
in $M/\sigma.$
\endproclaim
Now let $X_1, X_2\subset\P^n$ be a pair of smooth curves
meeting transversely at a point $p$, and set
$$X=X_1\cup X_2, M=N^*_X\otimes k(p)$$
with $\sigma\in M$ the singular element. For any vector bundle $E$
on a smooth curve $C$, we denote by $HN.(E)$ the (increasing)
Harder-Narasimhan filtration of $E$ and, for any point $p\in C$,
by $HN.(E, p)$ the fibre of the latter at $p$, i.e. $HN.(E)\otimes
k(p).$ We refer to $HN_1(E)$ as the {\it{positive subsheaf}} of
$E$ and denote it by $E_+$, and to $HN_1(E,p)$ as the
{\it{positive subspace}} of $E(p)=E\otimes k(p).$
\proclaim{Definition 2.3} $X_1$ and $X_2$ are said to have
{\it{good interface}} (or to {\it{interface well}}) at $p$ if the
filtrations
$$HN_\bullet(N^*_X|_{X_i}, p)\subseteq M, i=1,2$$
are in relative general position.\endproclaim For instance, in the
above example we showed a transverse pair of lines $L_1, L_2$ in
$\P^n$ do not interface well (at their point of intersection).This
is closely related to the fact that $L_1\cup L_2$ is not almost
balanced. Indeed the following general remark is easy to prove
\proclaim{Lemma 2.4} If $X_1, X_2$ is a general pair of rational
curves in $\P^n$ meeting at $p$ and both are almost balanced, and
$X=X_1\cup X_2$ then $N^*_X|_{X_i}$ is almost balanced for $i=1,2$
and $X_1, X_2$ have good interface at $p$ iff  the fibres at $p$
of the positive subbundles of $N^*_{X_i}, i=1,2$ meet transversely
in $N^*_X(p).$\qed\endproclaim

In the case where the ambient space is $\P^n$ it will often be
convenient to work with the twisted bundles $N\otimes L\inv,
N^*\otimes L$, where $L=f^*\O_{\P^n}(1)$ (which might be called
the normalized normal and conormal bundles, and it will be
convenient to use the notation
$$\N=N\otimes L\inv, \N^*=N^*\otimes L.$$

\subheading{3. The osculatrix filtration} \ss For a smooth scheme
$X/B$, and a (say locally free) sheaf $L$ on $X$, we denote by
$P^m_{X/B}(L)$ or just $P^m(L)$ the sheaf of $m$th order principal
parts of $L$ (cf. \cite{ EGA}). This sheaf carries a natural
increasing filtration with the $i$th quotient being $P^i(L)$ and
the $i$th graded piece being $\Sym^i(\Omega_{X/B})\otimes L$. We
will denote the $i$th subsheaf in this filtration, i.e. the kernel
of the natural map
$$P^m(L)\to P^i(L)$$
by $P^{[m,i)}(L)$. \remark{Example 3.1} If $X=\P^n, L=\O(1),
V=H^0(L)$, the canonical map
$$V\otimes\O_X\to P^1(L)$$
is easily seen to be an isomorphism.\qed\endremark Now let $C\to
\P^n$ be a smooth curve and again let $V=H^0(\O_{\P^n}(1)),
L=\O_C(1)$. We have a natural map
$$\rho:V\otimes\O_C\to P^{n-1}_C(L).$$
Note that the 'expected degeneracy' of $\rho$ is in
codimension 2. We will say that $C$ is {\it{totally
unramified}} if $\rho$ is surjective.\par
Now suppose that $C$ is totally unramified and set, as usual
$$M=M_C=\Omega^1_{\P^n}(1)\otimes\O_C.$$
Then in light of Example 3.1 $M$ coincides with the kernel of the
natural evaluation map
$$V\otimes\O_C\to L$$ hence we get a natural map
$$\rho: M\to P^{[n,0)}(L), P=P_C,$$
inducing a surjection
$$M\to P^{[n-1,0)}(L).$$Set
$$M^i=\rho\inv(P^{[n,i)}(L)), i=0,...,n-1.$$
Then $(M^.)$ is a descending filtration of $M$ with
$$M^{i-1}/M^i=K^i\otimes L, 1\leq i\leq n-1,\tag 3.1$$
where $K=\omega_C, $ while $M^{n-1}$ is a subsheaf, usually
proper, of $K^n\otimes L$ Note that by (3.1),
$$M^{n-1}\otimes L^{n-1}\otimes K^{\binom{n}{2}}
=\det(M)=L\inv$$
and hence
$$M^{n-1}=L^{-n}\otimes K^{-\binom{n}{2}}.$$
Note that by definition $M^1$ coincides with the twisted
conormal sheaf $N^*(L)$, whose fibre at any point $p\in C$
is the set of linear forms vanishing on the embedded tangent
line
$$T_pC\subset\P^n.$$ Similarly the fibre of $M^i$ at
$p$ for $i\leq n-1$ coincides with the set of linear forms
vanishing on the $i$th tangent (or $(i-1)$st osculating)
 space $T^i_pC$, i.e.
the set of linear forms vanishing to order at least $i+1$ at $p$
(the latter may be taken as the definition of $T^i_pC$); by our
assumption of total non-ramification $T^i_pC$ has constant
dimension for all $p\in C, i\leq n-1$ and if $i<n-1$ then $T^i_pC$
osculates to $C$ to order exactly $i+1$. \remark{Example 3.2: The
Rational Normal Curve} Let
$$C=C_n\subset\P^n$$
be the rational normal curve. Thus
$$V=H^0(\O_C(n))=:V_n,$$
hence the map $$V\otimes\O_C\to P^n(\O_C(n))$$
is fibrewise injective ('a degree-$n$ polynomial is
determined by its $n-$th order Taylor expansion at any point')
hence this map is in fact an isomorphism.
Likewise, we may identify the fibre of $M^i$
at $p\in C$ with
$$H^0(\O_C(n-(i+1)p))=H^0(\O_C(n-i-1)),$$
an identification which, up to scalars, is {\it{ independent
of the point}} $p\in C$. It follows that
$$M^i\simeq V_{n-1-i}\otimes L_i$$
where $V_{n-1-i}=H^0(\O_C(n-i-1))$ is an $(n-i)$-dimensional
vector space and $L_i$ is some line bundle on $C$. Comparing
degrees via (3.1) we see that $L_i$ has degree $-i-1$, i.e.
$$M^i\simeq V_{n-1-i}\otimes\O_C(-i-1).\tag 3.2$$
Alternatively, and more directly, one may observe that
given $p\in C\simeq\P^1$ and a linear form $\l_p$
on $\P^1$ zero at $p$, elements of the fibre at $p$ of
$V_{n-1-i}\otimes\O_C(-i-1)$ may be represented uniquely
in the form $f.\l_p^{i+1}, f\in V_{n-1-i}$ which gives rise
to a natural inclusion
$$V_{n-1-i}\otimes\O_C(-i-1)\to V_n;$$
the image of this inclusion at $p$ clearly coincides with
the set of polynomials whose $i$th jet at $p$ is zero, thus
the image globally coincides with $M^i$, so
$M^i\simeq V_{n-1-i}\otimes\O_C(-i-1)$ and
we have
exact sequences
$$0\to V_{n-1-i}\otimes\O_C(-i-1)\to V_n\otimes\O_C\to
P^i(\O_C(n))\to 0.\tag 3.3$$ In particular for $i=1$ we deduce the
isomorphism
$$N^*_{C_n}\simeq V_{n-2}\otimes\O_C$$
Note that under the above identification the inclusion
$M^i\subset M^j$ is given at $p\in C$ by
$$f\mapsto f\l_p^{j-i}$$
where $\l_p$ is the unique, up to scalars, linear polynomial
vanishing at $p$. From this it follows easily that
for any $i>j,$ the inclusion
 $M^i\subset M^j$ is 'nondegenerate' in the sense
 that its image
is not contained in a 'flat' subbundle of the form
$$W\otimes\O_C(-j-1), W\subsetneq V^j.$$
\endremark
\proclaim{Corollary 3.3} A general rational curve of degree $d$ in
$\P^n$ is totally unramified.\endproclaim \demo{proof} We may
assume the curve $C$ is nondenerate, i.e. $d\geq n$. Then $C$ is
obtained as the projection of a rational normal curve in $\P^d$
corresponding to a general $(n+1)$-dimensional subspace
$$V\subset H^0(\O_C(d))=:V_d.$$ Recall the
isomorphism
$$V_d\otimes\O_C\simeq P^d(\O_C(d))$$
and the natural surjection
$$P^d(\O_C(d))\to P^{n-1}(\O_C(d)).$$
Then $C$ being totally unramified is equivalent to $V$ mapping
surjectively to $P^{n-1}(\O_C(d)).$ That this holds for a general
$V$ can be seen by a standard (and trivial) dimension
count.\qed\enddemo \example{Example 3.4} Let $C_d\subset\P^n$ be a
generic rational curve of degree $d\leq n.$ The $C_d$ is a
rational normal curve within the span $\P^d$ of $C_d$. As is well
known,
$$T_{\P^n}|_{\P^d}\simeq T_{\P^d}\oplus (n-d)\O(1)$$
therefore
$$N^*_{C_d}\simeq N^*_{C_d/\P^d}\oplus
(n-d)\O_{\P^d}(1)|_{C_d},$$ thus
$$N^*_{C_d}\simeq (n-d)\O(-d)\oplus (d-1)\O(-d-2).\tag 3.4$$
\endexample
A fact closely related to the exact sequence (3.3) for $i=1$ is
the following observation, which is probably well known, but of
which
we shall subsequently require a more precise from%
\proclaim{Lemma 3.5} For $n> 0$, there is a canonical isomorphism on $\P^1$%
$$P^1(\O(n))\simeq V_1\otimes\O(n-1)\tag 3.5$$
whose $H^0$ composes via the natural map
$$V_n\to H^0(P^1(\O(n)))$$
to yield the usual co-symmetrization map $$\phi:V_n\to V_1\otimes
V_{n-1},$$
$$\l_1\cdots\l_n \mapsto \sum\limits_i\l_i\otimes
\l_1\cdots\hat{\l_i}\cdots\l_n.$$\endproclaim%
\demo{proof} In terms of homogeneous coordinates $X_0, X_1,$, the
map $\phi$ is given by%
$$\phi(f)=X_0\otimes\partial f/\partial X_0 +X_1\otimes\partial f/\partial
X_1.$$ Being a 1st-order differential operator in $f$, this
expression clearly descends to a canonical $\O_{\P^1}-$linear,
PGL-equivariant map $$P^1(\O(n))\to V_1\otimes\O(n-1).$$ To show
this map is an isomorphism it suffices to check it is surjective
at one point, e.g. $[1,0]$, where this is obvious.\qed\enddemo
\remark{Remark 3.6} The fact that $P^1(\O(n))\simeq 2\O(n-1)$ also
follows from the fact that the extension class of the natural
exact sequence
$$0\to\Omega_{\P^1}(n)\to P^1(\O(n))\to\O(n)\to 0$$
represents the first Chern class of $\O(n)$, hence the sequence is
nonsplit if $n\neq 0.$ However, the explicit isomorphism given
above will be important in the sequel.
\endremark
 \subheading{4. Realizing the splitting geometrically}\ss
Continuing with the case of the rational normal curve
$C=C_n\subset\P^n$, consider the special case of (3.2) that is the
isomorphism
$$N^*(L)\simeq V_{n-2}\otimes\O(-2)\tag 4.1$$
where $L=\O_C(n)=\O_{\P^n}(1)|_C$ and
$$N^*=M^1_C=\I_C/\I^2_C$$ is the
conormal bundle. This splitting of the conormal bundle
may be realized geometrically as follows.
Note that the set of divisors of degree $n-2$ on
$C$ may be identified with $\P(V_{n-2})$.
Let $\D$ be such a divisor.
Projection from $\D$, that is, from
its linear span $\bar{\D}\simeq\P^{n-3}\subset\P^n,$
maps $C$ to a smooth conic in $\P^2$, whose pullback
$K_\D$, i.e. the cone on $C$ with vertex $\bar{\D}$,
 is a rank-3 quadric in $\P^n$ containing $C$.
 Note that the map
 $$\D\mapsto K_\D$$
 is one-to-one because $\bar{\D}$ coincides with the
 singular locus of $K_\D$ while
 $$\D=\bar{\D}\cap C.$$
The equation of $K_\D$ gives rise to a nonzero section of
$N^*(2L)=N^*(2n)$ vanishing on $\D$. Untwisting by $\D$, we get a
subsheaf
$$\kappa_\D\simeq\O_C(-2)\subset N^*(L).$$
As $N^*(n+2)$ is a trivial bundle,
$\kappa_\D$ must be a (saturated) subbundle.
Now the natural map
$$H^0(\I_C(2L))\to H^0(N^*_C(2L))$$
is clearly injective because no quadric can be
double along $C$; indeed the singular locus of any
quadric must be a proper linear subspace of $\P^n$,
hence cannot contain $C$.
Thus the assignment
$$\D\mapsto\kappa_\D$$
is one-to one, giving rise to a one-to-one map
$$\lambda:\P(H^0(\O_C(n-2)))\to \P(N^*_C(n+2)).$$
Since both source and target of $\lambda$ are $\P^{n-2}$'s,
$\lambda$ must be a projective isomorphism, arising
from a linear isomorphism
$$\lambda:H^0(\O_C(n-2))\to H^0(N^*_C(n+2)).$$
It would be nice to contruct $\lambda$ directly as
a linear map of vector spaces, but we don't
know how to do it.
Moreover it is clear by construction that for a general
$p\in C$, the fibre $\kappa_\D(p)$ corresponds to the
hyperplane in $N_C(p)$ that comes from the hyperplane
$$<T_pC,\D>\subset\P^n$$
and it follows easily from this that for a general
choice of divisors $\D_1,...,\D_{n-1}$ and a general
$p\in C$, the subsheaves $\kappa_{\D_1},...,
\kappa_{\D_{n-1}}$ are independent at $p$, whence a
generically injective map
$$\phi:(n-1)\O(-2)\to N^*(L)$$
which, in view of (4.1), must be an isomorphism, giving the
desired realization of the splitting (4.1). Note that it follows a
posteriori that $\kappa_\D\subset N^*(L)$ is a subbundle
isomorphic to $\O_C(-2)$ for {\it{all}} divisors $\D$, not just
general ones. Note also that for any $\D$ we may choose
homogeneous coordinates $X_0,...,X_n$ so that $\kappa_\D$ is given
off $\D$ by
$$F=X_0X_2-X_1^2$$
or in terms of affine coordinates $x_1,...,x_n$, by
$$f=x_2-x_1^2.$$

\subheading {5. Vector bundles on some rational trees}
\medskip
In view of Grothendieck's theorem about decomposability of vector
bundles on nonsingular rational curves, it is natural to ask to
what extent decomposability holds for vector bundles on
{\it{rational trees}}, i.e. nodal curves of arithmetic genus 0.
The following result is hardly surprising. \proclaim{Proposition
5.1} Let $C$ be a nodal curve of the form
$$C=C_1\cup_pC_2$$
where $C_1,C_2$ are nonsingular, rational, and meet only at $p$
(we call such a curve $C$ a rational angle). Then any locally free
coherent sheaf on $C$ is a direct sum of line
bundles.\endproclaim\demo{proof} We begin with the following
observation. Let $F$ be a vector bundle on a nonsingular rational
curve $D$, $q$ a point on $D$ and $v\in F(q)$ a nonzero element of
the fibre at $q$. Then there is a basis of $F(q)$ that contains
$v$ and is compatible with the filtration induced on $F(q)$ by the
Harder-Narasimhan filtration $HN.(F)$. This is a triviality. It
follows from it that there is a line subbundle $L\subset F$ so
that $v\in L(q)$ and that
$$F\simeq L\oplus (F/L).$$ Indeed if $i$ is such that
$$v\in HN_i(F)(q)\setminus HN_{i-1}(F),$$ there is a unique line
bundle summand $\bar{L}$ of $HN_i(F)/HN_{i-1}(F)$ such that
$$v\in
\bar{L}(q)\mod HN_{i-1}(F),$$ and we can take for $L$ any lifting
of $\bar{L}$ to $$HN_i(F)\simeq HN_{i-1}(F)\oplus
HN_i(F)/HN_{i-1}(F)\subseteq F.$$ Now let $E$ be a vector bundle
on $C$ as in the Proposition, and let $E_1\subseteq E|_{C_1}$ be
the positive subbundle, i.e. the smallest nonzero subsheaf in the
Harder-Narasimhan filtration of $E|_{C_1}$.
 Pick any
$$0\neq v\in E_1(p)$$
and apply the above observation with $D=C_2, q=p, F=E|_{C_2}.$ It
yields a line bundle summand $L\subset F$ with fibre at $p$
generated by $v$ and a complementary summand $G\subset F$. Set
$$W=G(p)\subset E(p).$$
Now there is a line bundle summand $M$ of $ E_1|_{C_1}$ with fibre
at $p$ generated by $v$, and clearly $M$ is also a summand of
$E|_{C_1}$. Thanks to the fact that $v\not\in W$, it follows that
$W$ is in general position with respect to the filtration on
$E(p)$ induced by the Harder- Narasimhan filtration of $E|_{C_1}$.
Hence there is a complementary subbundle $B$ to $M$ in $E|_{C_1}$
with $B(p)=W.$ Then $L,M$ glue to a line subbundle $\Lambda$ of
$E$ and $G,B$ glue to a complementary subbundle $\Gamma$, with
$$E\simeq \Lambda\oplus \Gamma.$$ By induction on the rank of $E$,
$\Gamma$ is a direct sum of lines bundles, hence so is $E$.\qed

\enddemo
Next, recall that a vector bundle $F$ on a nonsingular rational
curve $D$ is said to be {\it{almost balanced}} if it has the form
$$F=r_+\O(k)\oplus s\O(k-1)$$
and {\it{balanced}} if we may moreover assume $s=0$. The subsheaf
$$F_+=r_+\O(k)\subseteq F$$
is then uniquely determined and called the {\it{positive
subsheaf}} of $F$. The following remark, whose (trivial) proof we
omit, gives a useful cohomological characterization of almost
balanced bundles. \proclaim{Lemma 5.2} F is almost balanced iff
some twist $G$ of $F$ satisfies
$$H^1(G)=0,\ \ H^0(G(-1)=0.\qed $$\endproclaim
A vector bundle $E$ on a rational angle $C=C_1\cup_pC_2$ is said
to be {\it{almost balanced}} or AB if it has  either the form
$$E=r_+\O(k,\l)\oplus s_1\O(k,\l-1)\oplus s_2\O(k-1,\l)
\tag 5.1$$ or the form
$$E=r_{1+}\O(k+1,\l)\oplus r_{2+}\O(k,\l+1)\oplus s\O(k,\l)
\tag 5.2$$ where $\O(a,b)$ denotes the lines bundle having degree
$a$ on $C_1,$ $b$ on $C_2$ (any line bundle on $C$ is one of these
for a unique $(a,b)$); $E$ is {\it{balanced}} if we may further
assume $r_+=0$ in (5.1) or $s=0$ in (5.2). A convenient
characterization of AB bundles is the following
\proclaim{Proposition 5.3} Given a vector bundle $E$ on a rational
angle, the following are equivalent:\par (i) $E$ is almost
balanced;\par (ii)the bundles $E_i=E|{C_i}, i=1,2$ are almost
balanced and the positive subspaces
$$(E_i)_+(p)\subseteq E(p), i=1,2$$
are in general position;\par (iii) some twist $G=E(-k,-\l)$
satisfies  either
$$H^1(G)=H^0(G(-1,0))=H^0(G(0,-1))=0\tag 5.3$$
or
$$H^0(G)=H^1(G(1,0))=H^1(G(0,1))=0\tag 5.3'$$
\endproclaim\demo{proof}
The fact that (i) implies (ii) and (iii) is trivial. The proof
that (ii) implies (i) is very similar to that of Proposition 5.1
(whether (5.3) or (5.3') occurs depends on whether the two
positive subspaces $(E_1)_+(p), (E_2)_+(p)$ span $E(p)$ or not).
To prove that (iii)- say in the form (5.3)-
 implies (i) we use
Proposition 5.1. We may assume $E$ contains $\O(0,0)$ as a direct
summand, in which case the $H^0$ vanishing hypotheses in (5.3)
show that $E$ cannot have a direct summand $\O(a,b)$ with either
$a>0$ or $b>0,$  while the $H^1$ vanishing hypothesis implies that
$E$ cannot have a direct summand $\O(a,b)$ with either $a<-1$ or
$b<-1$ or $(a,b)=(-1,-1)$. Hence $E$ is a sum of copies of
$\O(0,0), \O(0,-1)$ and $\O(-1,0)$ so it is almost balanced.
\qed\enddemo \proclaim{Corollary 5.4} If $(C_t,E_t)$ is a flat
family of pairs (curve, vector bundle) such that $C_0$ is a
rational angle, $E_0$ is almost balanced and a general $C_t$ is
smooth. Then a general $E_t$ is almost balanced.\endproclaim
\remark{Remark 5.4.1} A family of pairs (curve, vector bundle)
coming from a flat family $\Cal C/T$ and a vector bundle $\Cal E$
on $\Cal C$ is said to be {\it{quasi-constant}} if there is a
filtration of $\Cal E$ with vector bundle quotients which
restricts to the HN filtration on each fibre. Then the same
argument shows that if $E_0$ is almost balanced then any family is
quasi-constant: in fact, locally over the base, $\Cal E$ itself
splits as a direct sum of line bundles.\endremark
\medskip

\remark{Example 5.5: Piecewise degenerate rational angles} The
argument in the proof of Proposition 5.1 may be used to construct
splittings of vector bundles on rational angles. As an example,
which will be needed in the sequel, we consider a general rational
angle of the form
$$C=C_{a,b}=C_a\cup C_b\subset \P^n$$
where $C_a, C_b$ have degrees $a,b$ respectively
with $a,b <n$ and they meet at $p$.
 Thus $C_a$ is a rational normal curve in
$\P^a=<C_a>$ and likewise for $C_b$.
Set
$$N^*=N^*_C.$$
Suppose first
that
$$a+b<n.$$ Then $\P^a\cup\P^b$ spans a $\P^{a+b}$ and we
have a splitting
$$N^*\simeq N^*_{C/\P^{a+b}}\oplus (n-a-b)L^*,\
L^*=\O_{\P^n}(1)|_C=\O((-a)\cup(-b)).\tag5.4$$ Therefore it
suffices to analyze the case $a+b\geq n.$ Now clearly we have
$$N^*_{C_a}\simeq (n-a)\O(-a)\oplus (a-1)\O(-a-2)$$
and $N^*|_{C_a}$ is obtained from this by an elementary
modification at $p$ corresponding to $T_pC_b$,
which is a general direction. This modification has the
effect of chopping an $\O(-a)$ summand down to an $\O(-a-1),$
and it follows easily that
$$N^*|_{C_a}\simeq (n-a-1)\O(-a)\oplus
\O(-a-1)\oplus (a-1)\O(-a-2).$$
The HN filtration of this sheaf is
$$HN_1(N^*|_{C_a})=(n-a-1)\O(-a),$$
$$HN_2(N^*|_{C_a})=(n-a-2)\O(-a)\oplus\O(-a-1).$$
Ditto for $C_b.$
Let $M=N^*(p)$  and let
$$\Phi_i^c\subset M, i=1,2, c=a,b$$
be the fibres of these HN sheaves. Then $\Phi^a_1$
may be identified with the set of linear forms vanishing on
$\P^a\cup T_pC_b$. In particular it does not contain the
singular element $\sigma$ and its image in the Zariski
conormal space
$$S\simeq M/<\sigma>$$
is a generic $(n-a-1)-$dimensional subspace and ditto for
$b$. Since
$$(n-a-1)+(n-b-1)<n-1=\dim(S)$$
by our assumption $a+b\geq n$, these subspaces have zero
intersection
and in particular $\Phi^a_1\cap\Phi^1_b=0.$ On the other
hand, clearly
$$\Phi^2_a=\Phi^1_a\oplus <\sigma>$$
and likewise for $b$. Given this information, it is easy to see,
arguing as in the proof of Proposition 5.1, that we have the
following splitting (valid for $a+b\geq n$)
$$N^*\simeq (n-a-1)\O((-a)\cup (-b-2)\oplus
(n-b-1)\O((-a-2)\cup(-b))$$$$\oplus\O((-a-1)\cup (-b-1))
\oplus(a+b-n)\O((-a-2)\cup(-b-2)).\tag 5.5 $$ Note that if $a+b>n$
this bundle is never almost balanced. When $a+b<n$ we may apply
(5.5) to the inclusion $C_a\cup C_b\subset\P^{a+b}$ and in light
of (5.4)
 we get the following splitting in that case
$$N^*\simeq (n-a-b)\O((-a)\cup(-b))\oplus
(b-1)\O((-a)\cup (-b-2)$$$$\oplus
(a-1)\O((-a-2)\cup(-b))\oplus\O((-a-1)\cup (-b-1)) \tag 5.6$$

\endremark\ss

Despite Proposition 5.1, it is not in general true that a vector
bundle, even of rank 2, on a rational tree is decomposable, as the
following example shows. \remark{Example 5.6} Consider a nodal
curve of the form
$$C=C_0\cup C_1\cup\cdots C_4$$
where each component $C_i$ is a $\P^1$, $C_0$ meets each $C_i,
i>0$, in a unique point $p_i$ and there are no other
intersections. A vector bundle $E$ on $C$ may be constructed by
taking a copy of $\O(1)\oplus\O$ on each component $C_i$ and
gluing together generically at the $p_i$. Note that if a line
subbundle $L$ of $E$ has degree 1 on some $C_i, i>0$ then it has
degree 0 on $C_0$ and there is at most one other component $C_j$
on which $L$ has degree 1. Similarly, if $L$ has degree 1 on $C_0$
then it has degree 0 on every other $C_i$. It follows that any
line subbundle of $E$ must have degree at most 2 and since $E$ has
degree 5 it is indecomposable.\qed\endremark

In general, when an almost balanced bundle specializes to a
non-almost balanced one, there is no well-defined limit to the
maximal subbundle. In the next result, however, we identify one
very special case when the limit can be at least partly
identified. \proclaim{Proposition 5.7} Let $X/B$ be a proper
family with general fibre $\P^1$ and special fibre $X_0$ either a
$\P^1$ or a rational angle,  with $X$ a smooth surface and and $B$
a smooth curve. Let $E$ be a vector bundle on $X$ whose
restriction $E_b$ on a general fibre $X_b$ is almost balanced and
whose restriction $E_0$ on $X_0$ admits a filtration
$$E_{00}=0\subseteq
E_{01}\subseteq E_{02}\subseteq E_{03}=E_0$$
such that for some integer $k$, each $E_{0(i+1)}/E_{0i}$
splits as a direct sum of line bundles of total
degree $k-i$.Then\par
(i) if $\rk (E_{01})>\rk (E_{03}/E_{02}),$ the maximal
subbundle $E_{b+}$ specializes to a direct summand
 of $E_{0}$
and of $E_{01}$ that is a direct sum of line bundles of
total degree $k$;\par
(ii)if $\rk (E_{01})<\rk (E_{0}/E_{02}),$ the minimal
quotient $E_{b-}$ specializes to a direct summand
 of $E_0/E_{02}$ and of $E_0$
that is a quotient
bundle that is a direct sum of line bundles of
degree $k-2$;\par
(iii) if $\rk (E_{01})=\rk (E_{0}/E_{02}),$
then $E_b$ is balanced.\endproclaim\demo{proof}
Assertion (iii) is trivial and (i) and (ii)
are mutually dual, so
it will suffice to prove (i). To that end, note to
begin with that by simple arithmetic, if we set
$$r=\rk (E_{01})-\rk (E_{03}/E_{02})$$ then
the maximal subbundle
$$E_{b+}=r\O(k).$$
By unicity of the maximal subbundle $E_{b+}$, there exists
a subsheaf of $E$ which restricts on the generic fibre
$X_b$ to $E_{b+}$; let $E_+$ be the saturation of such
a subsheaf. Thus $E/E_+$ is torsion-free, hence by elementary
depth considerations, $E_+$ is locally free, hence its
restriction $E_{+0}$ on $X_0$ is a direct sum of $r$ line
bundles, whose total degrees add up to $rk.$\par
Now let $S\subseteq E_{0+}$ be a a line subbundle of
maximal total degree. This degree is clearly at
least $k$. Then $S^*\otimes E_0$ is a direct sum of line
bundles of nonpositive total degree and admits a regular
(locally nonzerodivisor) section. This clearly implies
$S$ is a direct summand of $E_0$ (and of $E_{01}$)
of degree exactly $k$. Therefore $E_{0+}$ is a direct sum of
line bundles of total degree exactly $k$ and
is a direct summand of $E_0$ and of $E_{01}.$
\qed
\enddemo
\subheading{6. Normal bundles of generic rational curves and
angles}\ss A smooth rational curve or rational angle
$C\subset\P^n$ is said to be {\it{almost balanced}} if its normal
bundle is. The {\it{bidegree}} of a rational angle $C_a\cup
C_b\subset\P^n$ is defined to be $(a,b)$. Our main purpose in this
section is to prove the following result \proclaim{Theorem 6.1}
\item{(i)} A generic rational curve of degree $d\geq n$ in $\P^n$
is almost balanced. \item{(ii)} Let
$$C=C_a\cup C_b\subset\P^n$$ be a generic rational angle
of given bidegree $(a,b)$.
Then $C_a, C_b$ interface well provided
$$\max(a,b)\geq n. $$\endproclaim
Before turning to the proof of Theorem 6.1 we note some explicit
corollaries. First some notation. Fixing $n$, define integers
$k(d), r(d)$ by
$$(n+1)d-2=k(d)(n-1)+n-1-r(d), 0<r(d)\leq n-1.$$
Note that an almost balanced bundle of degree $-(n+1)d+2$ must
have splitting type $((-k)^r,(-k-1)^{n-1-r})$, therefore
\proclaim{Corollary 6.2} A generic rational curve $C_d$ of degree
$d\geq n$ in $\P^n$ has conormal bundle
$$N^*=r(d)\O(-k(d))\oplus(n-1-r(d))\O(-k(d)-1).$$
\endproclaim

We can similarly determine the splitting type of the normal bundle
of generic rational angles. \proclaim{Corollary 6.3} For a generic
rational angle $C$ of bidegree $(a,b)$ in $\P^n$ with $a,b\geq n$,
$C$ is almost balanced and we have: if $r(a)+r(b)\geq n+1$,
$$N^*_C\simeq (r(a)+r(b)-n-1)\O((-k(a))\cup (-k(b)))\oplus$$
$$(n-r(b))\O((-k(a))\cup (-k(b)-1))\oplus(n-r(a))
\O((-k(a)-1)\cup (-k(b)));\tag 6.1$$ if $r(a)+r(b)\leq n,$
$$N^*_C\simeq (r(a)-1)\O((-k(a))\cup (-k(b)-1))\oplus
(r(b)-1)\O((-k(a)-1)\cup(-k(b)))$$
$$\oplus (n+1-r(a)-r(b))\O((-k(a)-1)\cup (-k(b)-1)).\tag 6.2$$
\endproclaim
\demo{proof} From almost balancedness of $C_a, C_b$ it
follows easily that, setting $N^*=N^*_C,$ we have
$$N^*|_{C_a}\simeq (r(a)-1)\O(-k(a))\oplus (n-r(a))
\O(-k(a)-1)$$ and likewise for $b$. Then almost balancedness of
$C$ is equivalent to the positive subspaces of these bundles at
the node $p$ being in general position. Since this holds,
analyzing $N^*$ as in the proof of Proposition 5.1 and Example 5.5
yields the claimed splitting.\qed\enddemo \proclaim{Corollary 6.4}
For a generic rational angle $C$ of bidegree $(a,b)$ in $\P^n$
with $1\leq a\leq n-1$, $n\leq b$, we have, setting $r=r(b),
k=k(b)$:\newline Case 1: if $r(b)>a$,
$$N^*_C\simeq (r(b)-a-1)\O((-a)\cup(-k))\oplus \O((-a-1)\cup
(-k))$$
$$\oplus(n-r(b))
\O((-a)\cup(-k-1))\oplus(a-1)\O((-a-2)\cup (-k)); \tag 6.3$$ Case
2: if $r(b)\leq a$
$$N^*_C\simeq (n-a-1)\O((-a)\cup(-k-1)\oplus\O((-a-1)\cup(-k-1))
$$$$\oplus(r(b)-1)\O((-a-2)\cup(-k))\oplus(a-r(b))
\O((-a-2)\cup(-k-1)).\tag 6.4$$ In particular, if $a=n-1$ then $C$
is almost balanced.
\endproclaim
\demo{proof} Analogous to the preceding proof, again using the
good interface of $C_a$ and $C_b$. Note that in the present case
$C$ is not necessarily almost balanced unless $a=n-1$, so we
cannot directly conclude from this that a general $C_{a+b}$ is
almost balanced.\qed\enddemo

Note that putting together Theorem 6.1, Corollaries 6.2-6.4 and
Example 5.5 we now know the splitting type of the normal bundle of
a generic rational curve or rational angle of every degree and
bidegree.\par
 The proof of Theorem 6.1 is somewhat long, so we break it into
steps. \subheading{Step 1: Case (n,1)} We begin by showing that
$C$ is almost balanced if $a=n,$ i.e. $C_a$ is a rational normal
curve, and $b=1$, i.e. $C_b$ is a line $L$ meeting $C_n.$ Applying
a suitable projective transformation, we may assume that $C_n$ is
the standard rational normal curve given parametrically in affine
coordinates by
$$x_i=t^i, i=1,...,n,$$
and that $L$ is a general line through the origin $p$. Note that
if we project parallel to the coordinates $x_4,...,x_n$ and prove
almost balancedness for the projected curve, it will imply almost
balancedness for the original; this is fairly clear a  priori, and
will become more clear with the computations that follow. Thus, it
suffices to prove $C$ is almost balanced in case $n=3.$ By
semi-continuity, it would suffice to prove $C$ is almost balanced
for one choice of $L$, and we pick the line with equation
$$x_1-x_3=x_2=0.$$
Now recall the identification (see \S3)
$$N^*_{C_3}=V_1\otimes\O(-5)$$
As usual, the restriction $N^*_C|_{C_3}$
is given by the elementary modification of
$N^*_{C_3}$ corresponding to the Zariski conormal
space to $C$ at $p$, and this conormal space
 is clearly generated by
the class of $x_2.$ Thus the positive subsheaf
$HN_1(N^*_C|_{C_3})$ is the unique 'special' (i.e. in this case,
degree- (-5)) subsheaf whose fibre at $p$ is generated by $x_2\mod
\m_{p, C_3}$, and clearly that subsheaf is generated locally by
$$f=x_2-x_1^2,$$ which is none other
than $\kappa_\D$ where $\D$ is the unique point at infinity
on $C_3,$ i.e. the point with homogeneous
coordinates $[0,0,0,1]$. Since a local basis for
$N^*_{C_3}$ is given by
$$f,g=x_2-x_1x_2$$ and $x_1$ is
a local parameter on $C_3,$ it also follows that a local
basis for the elementary modification $N^*_C|_{C_3}$
is given by $f, x_1g.$\par
Now over on the $L$ side, it is easy to see that
$N^*_C|_L$ has local basis $x_2,x_1(x_1-x_3)$ with
positive subsheaf $HN_1(N^*_C|_L)$ generated by $x_2.$
Almost balancedness for $C$ means that the two positive
subsheaves have different images in
the fibre $N^*_C\otimes k(p)$,
which amounts to saying that we modify $f$ to another
local section $f'$ of $N^*_C|_{C_3}$ that has the same
fibre at $p$, and that lifts to a local function vanishing
on $L$, viewing $f'$ as a local section of $N^*_C|_L$
and expressing it as a linear combination of
$x_2,x_1(x_1-x_3),$ the coefficient of $x_1(x_1-x_3)$
is nonzero at $p$.\par
 Indeed, set
$$f'=f+x_3g.$$
Since
$$x_3\in \m_{p,C_3}^3,$$
clearly $f$ and $f'$ have the same image in
$N^*_C|_{C_3}\otimes k(p).$ On the other hand, we have
$$f'=(1-x_1x_3)x_2-(x_1+x_3)(x_1-x_3)$$
and as section of $N^*_C|_L$, we have
$$f'=(1-x_1^2)x_2-2x_1(x_1-x_3),$$
which proves our assertion.\par \subheading{Step 2: More on (n,1)}
Note that for $n=3$ what we have proven, in fact, is that
$$N^*_C\simeq \O(-5\cup -2)\oplus\O(-6\cup -1)$$
and in particular $N^*_C$ is balanced. For $n>3,$ the positive
subsheaf $P$ of $N^*_C$ has corank 2 (i.e. rank $(n-3)$) and it is
at this point easy- and worthwhile too- to identify explicitly its
restriction on $C_n.$ Now in terms of the identification (4.1), we
have
$$HN_1(N^*_C|_{C_n}) = U\otimes\O(-n-2),$$
for some codimension-1 subspace $U\subset V_{n-2}$.
Clearly the intersection of $U$ with the open set of
divisors $\D$ not containing $p$ coincides with the
set of such divisors such that
$$L\subset\overline{\D+2p},\tag 6.6$$
because the quadric $K_\D$ is nonsingular at $p$ with tangent
hyperplane $\overline{\D+2p}$. Since both $U$ and the set of $\D$
satisfying (6.6) are linear spaces, it follows that, in fact
$$U=\{\D : L\subset\overline{\D+2p}\}.\tag 6.7$$
Now clearly
$$P|_{C_n}=W\otimes\O(-n-2)$$ where $W\subset U$
is a codimension-1 subspace (codimension 2
in $V_{n-2}$). Checking
again on divisors $\D$ not containing $p$, note
that if $$L\subset\overline{\D+p}$$ then the quadric $K_\D$
is smooth at $p$ and contains $L$. Its equation at $p$,
say $f$,
on the one hand clearly yields a section of
$HN_1(N^*_C|_{C_n})$; on the other hand, since the
tangent hyperplane to $K_\D$ at $p$ contains
$L\cup T_pC_n$, it is also clear that $f$
 yields a section of $HN_1(N^*_C|_L)$. Thus
$f$ in fact yields a section of $P$, so that
$$\D\in W.$$
Since the set of $\D$ with
$$p\not\in\D, L\subset\overline{\D+p}$$
is itself an open set in a codimension-2
subspace of $V_{n-2}$ (or its projectivization)
it follows that
$$\{\D:p\notin\D, \D\in W\}=\{\D:p\notin\D,
L\subset\overline{\D+p}\}$$
hence, finally, that
$$W=\{\D:L\subset\overline{\D+p}\}.\tag 6.8$$
The identity (6.8) has important consequences, including some
general position or genericity properties. For convenience, let's
temporarily set
$$M=N^*_C\otimes k(p), S:=<T_pC_n,L>^\perp=T_pC^\perp$$
where $T_p$ denotes the embedded
Zariski tangent space,
and note that $M$ contains 2 canonical hyperplanes, viz.
$$U=HN_1(N^*_C|_{C_n}, p)$$
and its analogue from the $L$ side,
$$U'=HN_1(N^*_C|_{L}, p),$$
as well as the 1-dimensional 'singular' subspace
$\sigma$ that is the kernel of either of the natural maps
$$M\to N^*_{C_n}\otimes k(p), M\to N^*_L\otimes k(p).$$
We have shown that
$$U\neq U'$$
and of course
$$\sigma\not\in U\cup U'.$$
Now note that $N^*_C|_L$, being an elementary
modification of $N^*_L$, and hence also $M$,
depend only on the flag
$$(p,L, S).$$
Consequently the group $G$ of projective motions preserving
this flag acts on $M$, preserving $U', \sigma$. Moreover it
is easy to see that $G$ acts transitively on the set of
hyperplanes in $M$ different from $U'$ and not containing
$\sigma$. The upshot is that for given flag $(p,L,S)$, $U$
may be assumed to be a generic hyperplane in $M$.
\par
Another important general position property that follows from
(6.8) is the following. Let $A\subset V_{n-2}$ be any irreducible
subvariety. Then the locus
$$\tilde{A}:=\{ (\D, p, L)\in A\times C_n\times \Bbb G(1,n):
p\in L\subset\overline{D+p}\}$$ maps to $A\times C_n$ so that
every fibre is a $\P^{n-3}$. Therefore $\tilde{A}$ is irreducible
and $(\dim(A)+n-2)-$dimensional. Therefore the fibre of
$\tilde{A}$ over a general pair $(p,L)$ with $p\in L$ is purely
$(\dim(A)-2)-$ dimensional. When $A$ is a linear space, so is this
fibre, and we conclude \proclaim{Corollary 6.5} Given any linear
subspace
 $A\subset V_{n-2}$, if $p\in C_n$ is sufficiently
general and $L$ is a sufficiently general line through
$p$, then $(p,L)$ impose independent conditions on
$A$ in the sense that
$$A(-(p,L)):=\{\D\in A:L\subset \overline{\D+p}\}$$
is a codimension-2 subspace of $A$ if $\dim(A)\geq 2$
and zero if $\dim(A)=1.$
\endproclaim
Translating this result into the language of normal bundles, we
conclude the following. \proclaim{Corollary 6.6} Assumptions as in
Corollary 6.5, if $\dim(A)\geq 2$ (resp. $\dim(A)=1$) then, the
intersection
$$A\otimes\O_{C_n}(-n-2)\cap
 N^*_{C_n\cup L}|_{C_n}$$
performed inside $ N^*_{C_n}$ and
viewed as subsheaf of  $ N^*_{C_n\cup L}|_{C_n}$,
extends to a subbundle
of $N^*_{C_n\cup L}$ of the form
$$A(-(p,L))\otimes(\O_{C_n}(-n-2)\cup\O_L(-1))
$$$$\oplus
\O_{C_n}(-n-2)\cup\O_L(-2)\oplus
\O_{C_n}(-n-3)\cup\O_L(-1)$$
(resp.
$$\O_{C_n}(-n-3)\cup\O_L(-2)).$$
\endproclaim

\subheading{Step 3: Case (n,1,1)}  Now our strategy for the proof
of the general case of Theorem 6.1 is essentially to degenerate
(implicitly) a curve with 1 or 2 components in $\P^n$ to a chain
consisting of lines, rational normal curves $C_n$ and some
degenerate rational normal curves $C_a$, each lying in an
$a-$plane in $\P^n$. To this end we consider next the case
$$C=C_n\cup L_1\cup L_2$$
where $L_1, L_2$ are generic lines among those incident to $C_n$
and meet it at generic points $p_1,p_2,$ respectively. We begin
with the case $n=3.$ In this case, note that
$$N^*_C|_{C_3}\simeq M\otimes\O(-6)$$
where $M$ is a 2-dimensional vector space canonically
isomorphic, up to scalars, to
$$M_i=N^*_C\otimes k(p_i)=N^*_{C_3\cup L_i}\otimes k(p_i),
i=1,2.$$
As above, $M_i$ contains a codimension-1 subspace $U'_i$
coming from the maximal subsheaf $\O(-1)$ on $L_i$ and,
by the general position property discussed above, we may
assume that under the identifications
$$\P(M)=\P(M_1)=\P(M_2),$$
we have
$$U'_1\neq U'_2.$$
This clearly implies that
$$N^*_C\simeq \O(-6\cup-1\cup-2)\oplus \O(-6\cup-2\cup-1).$$
Next consider the case $n=4$. Then with notations as above,
we have hyperplanes
$$U_i, U'_i\subset M_i, i=1,2$$
and a 1-dimensional subspace
$$W_i=U_i\cap U'_i.$$
On the other hand we clearly have
$$N^*_C|_{C_4}=\O(-6,-7,-7)$$
where the fibre $Z$ of the maximal subsheaf $\O(-6)$
coincides with $U_1\cap U_2$ under the natural embeddings
$$U_1, U_2\subset V_2.$$ Choosing things generally, we
may clearly assume
$$Z\not\subset U'_i, i=1,2.$$
Therefore $N^*_C$ admits a direct summand of the form
$$\O(-6\cup-2\cup-2).$$ Now let $a, b$ be a basis of $U_1'$
and lift them to subsheaves
$$A,B \simeq \O(-7)\subset N^*_C|_{C_4},$$
which is clearly possible. Modifying $A,B$ by the unique
up to scalar map
$$\O(-7)\to \O(-6)\subset N^*_C|_{C_4}$$
vanishing at $p_1$ we may assume
$$A\otimes k(p_2), B\otimes k(p_2)\subset U'_2.$$
Then $A,B$ project $\mod \O(-6)$
to generically, hence everywhere,
linearly independent $\O(-7)$ subsheaves, and
each of them glues to an
$\O(-1)$ subsheaves on $L_1$ and $L_2$.
Therefore finally
$$N^*_C\simeq \O(-6\cup-2\cup-2)\oplus 2\O(-7\cup-1\cup-1).
\tag 6.9$$\par Finally consider the case $n\geq 5.$ choosing $C_n,
L_1, L_2$ generically, we get as above hyperplanes
$$U_i\neq U'_i\subset M_i, i=1,2$$
where $U_i$ maps isomorphically to
$$\{\D:L_i\subset\overline{\D+2p_i}\}\subset V_{n-2},$$
under the natural map
$$N^*_C|_{C_n}\to N^*_{C_n},$$
and subspaces
$$W_i=U_i\cap U'_i,$$
which map isomorphically to
$$\{\D:L_i\subset\overline{\D+p_i}\}\subset V_{n-2}, i=1,2.$$
as well as a codimension-2 subspace
$$Z=U_1\cap U_2\subset V_{n-2}$$
such that
$$N^*_C|_{C_n}\simeq Z\otimes\O(-n-2)\oplus 2\O(-n-3).$$
Moreover
$$Y=Z\cap W_1\cap W_2$$
is of codimension 2 in $Z$ (hence vanishes if $n=5$).
Analyzing as above, we conclude
$$N^*_C\simeq Y\otimes\O(-n-2\cup-1\cup-1)
\oplus
$$$$ \O(-n-2\cup-1\cup-2)\oplus\O(-n-2\cup-2\cup-1)
\oplus 2\O(-n-3\cup-1\cup-1)\tag 6.10$$

\subheading{Step 4: The critical range}  Our next goal is to prove
that $C_d$ is almost balanced in the critical range $n<d\leq 2n.$
This range is difficult because when $d$ is in it, a rational
angle $C_{a,b}, a,b >1$ that is a limit of $C_d$ will usually not
be almost balanced. Still, for $d=n+1, n+2,$ note that almost
balancedness of $C_d$ follows by specialization from almost
balancedness of $C_n\cup L$ and $C_n\cup L_1\cup L_2$. In the
general case we will work inductively. By specialization again, it
will suffice to prove in the range $n+2\leq d\leq 2n-1$ that if
$C_d$ is balanced, then so is a general connected union $C_d\cup
L$. To this end, we will specialize $C_d$ to $C_n\cup C_{d-n}$ so
that $L$ specialized to a unisecant of $C_n$. Set
$$a=d-n\in [2,n-1], N^*=N^*_{C_n\cup C_{a}}$$
and note that
$$N^*|_{C_n}\simeq (n-2)\O(-n-2)\oplus\O(-n-3),$$
$$N^*|_{C_{a}}\simeq (n-a-1)\O(-a)\oplus\O(-a-1)
\oplus (a-1)\O(-a-2),$$ where the positive subsheaf
$(n-a-1)\O(-a)\subset N^*_{C_a}$ corresponds to the linear forms
vanishing on $C_a\cup T_pC_n$. As we have seen, choosing the span
of $C_a$ sufficiently generally, the fibre at $p$ of this positive
subsheaf does not contain the fibre at $p$ of the positive
subsheaf of $N^*|_{C_n}$, and consequently we have, when $a\leq
n-2,$ a decomposition
$$N^*\simeq(n-a-2)\O(-n-2\cup -a)\oplus \O(-n-3\cup-a)\tag 6.11
$$
$$\oplus \O(-n-2\cup-a-1)\oplus (a-1)\O(-n-2\cup-a-2).$$
If $a=n-1,$ we have a decomposition
$$N^*\simeq \O(-n-3\cup-n)\oplus(n-2)\O(-n-2\cup-n-1)$$
and in particular this bundle is (almost) balanced; this case is
similar to but simpler than the case $a\leq n-2$, so assume the
latter.

 Now let $P$ be the positive
subsheaf of $N^*_{C_d}$, which by almost balancedness is of the
form $\rho\O(k)$, and let $P'$ be the elementary modification of
$P$ at a general point $q\in C_d$ corresponding to a general line
$L$ through $q$, and $P_+\subset P'$ its positive subsheaf, which
is of the form $(\rho-1)\O(k)$ if $\rho >1$ and $(n-1)\O(k-1)$ if
$\rho=1.$ To prove that a $C_d\cup L$ is almost balanced ,
 it would suffice to prove that if $\rho >1$ then
the fibre at $q$ of $P_+$ is not contained in the fibre at $q$ of
the positive subsheaf of the elementary modification of $N^*_L$
corresponding to $T_qC_d$ (which corresponds to the set of
hyperplanes containing $L\cup T_qC_d,$ i.e. the Zariski conormal
space to $C_d\cup L$ at $q$). \par

Now let $P_0$ be the limit of $P$ on $C_n\cup C_a$, as computed in
Proposition 5.7. If $2a<n-1$, then the restriction of $P_0$ on
$C_n$ is of the form $\rho\O(-n-2)$ and our assertion follows from
Corollary 6.6. If $2a>n-1$, then $P_0|_{C_n}$ contains a subbundle
of the form $(2a-1)\O(-n-2)$ and since we may assume $a>1$ our
assertion again follows from Corollary 6.6. If $2a=n-1$ then
$P_0|_{C_n}$ contains $(n-2)\O(-n-2)$ so again we are done except
if $n=3$, in which case $a=1,$ against our assumption.
(Alternatively, almost balancedness of $C_{n+a}$ and $C_{n+a}\cup
L$ for $2a<n-1$ could also be established by degenerating
$C_{n+a}$ to $C_n$ plus $a$ general unisecants, but we shall not
need this.)\par Let us say that a bundle $E$ on a curve $C$ is
$b-$balanced, for an integer $b$, if it is a direct sum of line
bundles whose total degrees on $C$ is contained in an interval of
length $b$. The in the above proof we have shown that $C_n\cup
C_a$ is 2-balanced but not 1-balanced if $2\leq a\leq n-3,$
1-balanced if $a=1$ or $n-1$ and 0-balanced if $a=n-2$.

I claim next that
$$C=C_{n-1}\cup L\subset \P^n$$
is almost balanced. To see this, note that
$$N^*_C|_{C_{n-1}}\simeq \O(-n)\oplus(n-2)\O(-n-1),$$
$$N^*_C|_L\simeq (n-2)\O(-1)\oplus \O(-2)$$
Where the positive subsheaf $\O(-n)\subset N^*_C|_{C_{n-1}}$
is the chopped form of the subsheaf
$$\O(-n+1)\subset N^*_{C_{n-1}},$$
corresponding to
the unique hyperplane containing $C_{n-1}$. Thus, the fibre
at the node $p$ of this subsheaf is spanned by the singular
element, and that fibre
is not contained in the fibre at $p$ of
the positive subsheaf of $N^*_C|_{L}$. Therefore $C$
is almost balanced (even balanced) with normal bundle
$$N^*_C\simeq \O(-n\cup-2)\oplus (n-2)\O(-n-1\cup-1).$$
\subheading{Step 5: A gluing lemma} Our proof of almost
balancedness in higher degrees will be inductive, based largely on
the following result \proclaim{Lemma 6.7} Assume that
$$C_a\cup L_a, C_b\cup L_b\subset \P^n$$
are almost balanced, where $L_a, L_b$ denotes general
unisecant lines to $C_a$ or $C_b$ respectively. Then
$$C_a\cup L\cup C_b\subset\P^n$$
is almost balanced, where $L$ denotes a general line
meeting $C_a$  and $C_b$ .
\endproclaim\demo{proof} Set
$$C=C_a\cup L\cup C_b, N^*=N^*_C.$$
By almost balancedness of $C_a\cup L$, it is easy to see
that $N^*|_{C_a\cup L}$, which is a general elementary
modification of $N^*_{C_a\cup L}$ at a general point of
$L$, is of the form either
(case a0)
$$(\rho)\O(k\cup -1)\oplus 2\O(k\cup -2)\oplus (n-3-\rho)
\O(k-1\cup-1), \rho\geq 0,$$
or, if the positive subsheaf of $N^*|_{C_a}$ is of rank 1
(case a1),
$$\O(k\cup -2)\oplus (n-3)\O((k-1)\cup -1)\oplus
\O((k-1)\cup -2)$$
or, if  $N^*|_{C_a}$ coincides with its positive subsheaf,
i.e. is balanced
(case a2),
$$(n-3)\O(k\cup -1)\oplus 2\O(k\cup -2).$$
Likewise, we have analogous cases b0,b1,b2, where the integer
analogous to $k$ will be denoted by $\l$. Assume
first that neither a1 nor b1 hold.

Note that the positive subsheaf of $N^*|_{C_a\cup L}$,
restricted on $L$, meets the positive subsheaf
$(n-2)\O(-1)$ of $N^*|L$ in a subbundle of the
form
$$U_a\otimes\O(-1)=u_a\O(-1),$$ likewise for $b$. Applying
a suitable projective transformation to $C_b$, the subspaces
$$U_a, U_b\subset H^0(I_L(1))$$
clearly may be assumed in general position.
Then it is easy to see that if
$$u_a+u_b\geq n-1$$
then $N^*$ splits as
$$N^*\simeq (U_a\cap U_b)\otimes\O(k\cup-1\cup\l)\oplus$$$$
\nu_a\O(k\cup-1\cup(l-1))\oplus\nu_b\O((k-1)\cup -1\cup
\l)\oplus 2\O(k\cup-2\cup\l),$$
where
$$\nu_a=u_a-\dim U_a\cap U_b=n-1-u_b,$$
likewise for $\nu_b.$
Otherwise, $U_a\cap U_b=0$ and
$$N^*\simeq U_a\otimes\O(k\cup -1\cup(\l-1))
\oplus U_b\otimes \O((k-1)\cup-1\cup\l)$$$$
\oplus (n-3-u_a-u_b)\O((k-1\cup-1\cup(\l-1))
\oplus \O(k\cup-2\cup\l).$$
If case a1 holds but b1 does not, we have a splitting
$$N^*\simeq U_b\otimes\O((k-1\cup-1\cup\l)$$$$
\oplus((n-3-u_b)\O((k-1)\cup-1\cup(\l-1))\oplus
\O(k\cup-2\cup(\l-1))\oplus\O((k-1)\cup-2\cup\l).$$
Likewise if b1 but not a1 hold. Finally if a1 and b1 hold,
then
$$N^*\simeq(n-2)\O((k-1)\cup-1\cup(\l-1))\oplus
\O(k\cup-2\cup(\l-1))\oplus\O((k-1)\cup-2\cup\l).$$ This completes
the proof of Lemma 6.7.

\enddemo
\subheading{Step 6: Conclusion} As one consequence of Lemma 6.7,
we can now prove the almost balancedness of $C_d$ and $C_d\cup L$
in $\P^n$ for all $d\geq n$. The proof is by induction and the
case $d\leq 2n$ has been done previously. If $d>2n,$ specialize
$C_d$ to $C_n\cup M\cup C_{d-n}$ and $L$ to a general unisecant of
$C_n$. By Lemma 6.7 and induction, $C_n\cup M\cup C_{d-n}$ is
almost balanced. Moreover considering the known results about
$C_n\cup M\cup L$, it is easy to see that $L\cup C_n\cup M\cup
C_{d-n}$ is almost balanced as well. Therefore $C_d$ and $C_d\cup
L$ are almost balanced, completing the induction step.\par Next we
show that
$$C_a\cup C_b\subset\P^n$$ is almost balanced
whenever
$$n-1\leq a, n\leq b.$$
Indeed, it suffices to degenerate $C_a\cup C_b$ to $C_a\cup L\cup
C_{b-1}$ and use the fact that $C_a\cup L$ and $L\cup C_{b-1}$ are
almost balanced, together with Lemma 6.7.\par To complete the
proof of Theorem 6.1 it now suffices to show that $C_a,C_b$
interface well when
$$C=C_a\cup C_b,\ \ 2\leq a\leq n-2, b\geq n.$$
Then
we have
$$N^*|_{C_a}\simeq (n-1-a)\O(-a)\oplus \O(-a-1)
\oplus (a-1)\O(-a-2),$$
$$N^*|_{C_b}\simeq (r(b)-1)\O(-k)\oplus (n-r(b))\O(-k-1),
k=k(b).$$
If $r(b)=1$, good interface is automatic so assume $r(b)>1.$
As usual, set
$M=N^*\otimes k(p)$ and also let $ L=T_pC_a$ which, vis-a-vis
$C_b$, is a general unisecant at $p$.
Then $M$ splits naturally as
$$M=<\sigma>\oplus S$$
where $S$ may be canonically identified with the set of
linear froms vanishing on $L\cup T_pC_b.$
Now
$$S_1:=
HN_1(N^*|_{C_a},p)$$
may be identified with the set of
linear form vanishing on $C_a\cup T_pC_b,$ or equivalently
on  $\P^a\cup T_pC_b$
where $\P^a$ is the linear span of $C_a$,
hence for fixed $S$ and
$C_b$, $S_1$ may be considered a generic
$(n-1-a)$-dimensional
subspace of $S$. On the other hand, in proving $C_b\cup L$ is
almost balanced for a generic line
$L$ through $p$ we showed that
$$HN_1(N^*|_{C_b},p)\not\subset S.$$ Therefore
$$S_2:= S\cap HN_1(N^*|_{C_b},p)$$ is
$(r(b)-2)$-dimensional and meets $S_1$ transversely within $S$,
and therefore $S_1$ and $HN_1(N^*|_{C_b},p)$ meet transversely
within $M$.\par Next, $HN_2(N^*|_{C_a},p)$ is spanned by $S_1$ and
$\sigma$ and since $S$ maps isomorphically to $M/\sigma$ and $S_1$
is generic in $S$, it follows the the image of
$HN_2(N^*|_{C_a})\otimes k(p)$ in $M/\sigma$ is a generic
$(n-d-1)-$dimensional subspace, hence meets the image of
$HN_1(N^*|_{C_b})\otimes k(p)$ transversely. Thus $C_a$ and $C_b$
interface well at $p$.\qed

\subheading {7. Remarks on degenerating linear systems}\ss

The purpose of this section is to work out a convenient local
model for a specialization of a linear system as $\P^1$
specializes to a rational angle. Consider in $\P^1\times\P^1\times
\A^1$ the divisor $X$ (of type (1,1)) with equation
$$U_1V_1=sU_0V_0.\tag7.1$$

It is immediate that $X$ is a smooth surface and the projection
$$\pi:X\to B:=\A^1$$ is flat  with fibres $\pi\inv(s)=\P^1, s\neq 0$
and
$$X_0=\pi\inv(0)=X_1\cup X_2=(V_1=0)\cup (U_1=0)$$ a rational angle.
Moreover either
one of the projections $X\to\P^1\times\A^1$ is a blowup at one
point. Set
$$L=p_{\P^1\times\P^1}^*(\O(a,b)).$$
Then $L$ is a line bundle with degree $d=a+b$ on a fibre of $\pi$
and bidegree $(a,b)$ on $X_0.$ Moreover it is immediate from the
defining equation (7.1) that $\pi_*(L)$ is a trivial bundle of
rank $d+1$ with basis
$$U_0^aV_0^b, U_0^{a_0}U_1^{a-a_0}V_0^b,
U_0^aV_0^{b_0}V_1^{b-b_0}, a_0=0,...,a-1, b_0=0,...b-1,\tag7.2$$
or
in affine coordinates $u=U_1/U_0, v=V_1/V_0,$%
$$1,u,...,u^a, v,...,v^b.$$
 It
will be important for our purposes to determine the 'dualizing
principal parts' sheaf $P^+_{X/B}(L)$ (cf. (2.4)), which coincides
with the locally free double dual of
$P^1_{X/B}(L)$ :
\proclaim{Lemma 7.1} We have a linear isomorphism
$$\nabla^\l=(\nabla_1, \nabla_2):P^+_{X/B}(L)\overset\sim\to\to
\O_X(a-1,b)\oplus\O_X(a,b-1)$$
which, via the inclusion $P^1_{X/B}(L)\subset P^+_{X/B}(L)$,
corresponds to the differential operator $\nabla$ on $L$ which on
relative global sections is given by
$$\nabla^h:F\mapsto (\partial F/\partial U_0+sV_0/U_1\partial F/\partial V_1,
\partial F/\partial V_0+sU_0/V_1\partial F/\partial U_1)\tag7.3$$
$$=(\partial F/\partial U_0+V_1/U_0\partial F/\partial V_1,
\partial F/\partial V_0+U_1/V_0\partial F/\partial U_1)\tag7.4$$

and which in affine coordinates $(u,v)$ at $p$ is given by
$$\nabla^a:f\mapsto (af-u\partial f/\partial u+v\partial f/\partial v,bf+
u\partial f/\partial u-v\partial f/\partial v)\tag7.5
$$\endproclaim \demo{proof} To begin with, it is easy to check
that (7.3) and (7.4) are in fact equal, dehomogenize to (7.5) and
that they vanish on $F_0=U_1V_1-sU_0V_0$, hence send a multiple of
$F_0$ to another multiple of $F_0$. Therefore they define an
$\O_B-$linear, $\O_X(a-1,b)\oplus\O_X(a,b-1)$-valued differential
operator $\nabla$ on $L$ over $X$, a priori a rational one but
from the equality of (7.3) and (7.4) clearly $\nabla$ is regular.
Hence it defines a linear map $\nabla^\l$ on $P^1_{X/B}(L)$. This
map takes the value $(a,b)$ on $F=U_0^aV_0^b,$ i.e. $f=1$. Now
this section $f=1$ gives rise to a
trivialization of $L$ and hence to a splitting, locally near the node $p$,%
$$P^1_{X/B}(L)\simeq \O_X\oplus \Omega_{X/B},
P^+_{X/B}(L)\simeq \O_X\oplus \omega_{X/B}.$$ By 7.5, we have%
$$\nabla^\l(du)=(-u,u), \nabla^\l(dv)=(v,-v).$$

Now the dualizing sheaf $\omega_{X/B}$ is generated locally at $p$
by a form $\eta$ equal to $du/u$ on $X_1$ and to $-dv/v$ on $X_2$,
and clearly $$\nabla^\l(\eta)=(-1,1).$$ Since $(a,b)=\nabla^\l(1)$
and $ (-1,1)$ form a local basis for the target of $\nabla^\l,$ it
follows that $\nabla^\l$ extends to $P^+_{X/B}(L)$ to yield an
isomorphism as stated in the Lemma, locally near $p$. It is easy
to check that $\nabla^\l$ is an isomorphism locally at any
non-critical point, therefore it is an isomorphism.\qed

\enddemo

\subheading{8. Smoothing the normal bundle of a rational angle}\ss

In \S6 we determined the normal bundle to a generic rational angle
$$C_{a,b}\subset\P^n$$ and saw, in particular, that it is often not
almost balanced when $$\min(a,b)<n.$$ As a result, the locus of
these rational angles will appear as an improper part of the locus
of curves of degree $d=a+b$ with unbalanced normal bundle. For our
enumerative purposes, this locus must therefore be subtracted off,
{\it{with the correct multiplicity}}, to get the correct count of
{\it smooth} curves with unbalanced normal bundle. Our next goal,
then, is the determination of these multiplicities. This will
require a study of the deformation of the normal bundle in a
smoothing of $C_{a,b}$, which we will do using the methods of the
previous section.

Now fix natural numbers $n,a,b$ with $b\leq a, b<n$, set $d=a+b$
and consider a general local linear system of rank $(n+1)$ and
bidegree $(a,b)$, which we may consider extended to a subbundle
$W$ of $\pi_*(L)$ defined near 0. As a relative linear system on
$X/B$, $W$ is clearly very ample on $X_0$, with image a general
rational angle $C_{a,b}\subset\P^n$. Therefore $W$ is relatively
very ample in a neighborhood of $X_0$ and embeds a general fibre
$X_s$ as a smooth rational curve of degree $d$ in $\P^n$, thus
providing an explicit smoothing of $C_{a,b}$. We denote by
$$f=f_W:X\to\P^n$$ the associated mapping, defined near $X_0.$
Working in an affine neighborhood of the node $p=([1,0],[1,0])$ of
$X_0,$ we use coordinates $u=U_1/U_0, v=V_1/V_0.$\par

Now let $N$ denote the relative lci normal bundle for $f$ relative
to $\pi$ and $\N=N\otimes L\inv$ as usual. Fix a line bundle
$\tau$ of bidegree $(r+1, 1)$ on $X_0$ (note that $\tau$ is in
fact unique), and set
$$r=\lceil\frac{2d-2}{n-1}\rceil-1,$$
$$N_0=N\otimes L\inv(-\tau).$$
We say that the pair $(d,n)$ is {\it{perfect}} if $(n-1)|2(d-1)$.
The enumerative results of \S9 will apply only to perfect pairs.
The main result of this section is

\proclaim{Proposition 8.1}(Cohomological quasitransversality) if
$(d,n)$ is perfect, then for the generic smoothing $W$, $
R^1\pi_*(N_0)_0$ is killed by $\m_{0,B}$ and is a $\k(0)-$ vector
space of dimension $b-1$.

\endproclaim
\remark{Remark} The proof will show more generally, in the
not-necessarily-perfect case, that  $R^1\pi_*(N_0)_0$ is a direct
sum of $(b-1)\k(0)$ and a free part of rank
$(n-1)(r+1)-(2d-2)$.\endremark\demo{proof} To begin with, let us
dispense with the easy case where $N_0$ is almost balanced. In
that case, $N_0$ itself is a direct sum of line bundles $K$, and
the only ones which contribute to $R^1$ are $\O(-2,0)$ and
$\O(-1,-1)$, whose $R^1$ is easily seen to be free; moreover these
summands cannot occur at all in the perfect case.\par Turning now
to the proof proper, it will be based on the exact sequence (2.5)
which in our case takes the form
$$0\to N^*(L)\to W\otimes\O_X\to\O(a-1,b)\oplus\O(a,b-1)\to q\to
0\tag 8.1$$ with $q$ a (generic) length-1 skyscraper sheaf at $p$,
and which dualizes to
$$ 0\to\O(-a+1, -b)\oplus\O(-a, -b+1)\to W^*\otimes\O_X\to \N\to
q\to 0.\tag8.2$$ 
Twisting (8.2) by $\O(-r-1,-1)$ yields
$$0\to\O(-a-r, -b-1)\oplus\O(-a-r-1, -b)\overset\phi\to\to
W^*\otimes\O_X(-r-1, -1) \to N_0\to q\to 0.\tag8.3$$ Now by
Theorem 6.1 we have for the generic smoothing that
$$\pi_*(N_0)=0,$$
hence $R^1\pi_*(N_0)$ consists of a locally free part of rank
$(n-1)(r+1)-(2d-2)$ plus some torsion at 0. Moreover, it is easy
to check from the results of \S6 that
$$h^1(N_0|_{X_0})=((n-1)(r+1)-(2d-2))+(b-1)$$
or equivalently,
$$h^0(N_0|_{X_0})=b-1
$$ (note that $h^0(N_0|_{X_1})=0$). Therefore the value given
by Proposition 8.1 is, in an obvious sense, the smallest possible
value for $R^1\pi_*(N_0)$-- which explains the term 'cohomological
quasitransversality', and consequently by semi-continuity it will
suffice to exhibit one smoothing $W$ for which the assertion of
Proposition 8.1 holds.\par

The following Commutative Algebra assertion will allow us to
reduce to the case $r\leq 5$: \proclaim{Lemma 8.2} Let $A$ be a
regular local ring and $s,t$ part of of a regular system of
parameters on $A$. Let $M$ be an $A-$ module of finite type such
that for some natural number $n$, \par (i) $M/tM\simeq
nA/(t)$;\par (ii) $M/sM\simeq nA/(s,t).$\par Then $tM=0$
 or equivalently, $M\simeq
nA/(t).$\endproclaim\demo{proof} Use induction on $n$, which
coincides with the number of minimal generators of $M$. Suppose
first that $n=1,$ so $M=A/I$. By (i), we have $I\subseteq (t)$, so
write $I=tJ.$ Then by (ii), we have $ (tJ,s)=(t,s)$ so clearly
$J=(1)$ as desired. \par In the general case pick a primary cyclic
submodule $$N=Au\simeq A/Q\subseteq M,$$ where $u\in M$ is a
suitable minimal generator. Since $u$ maps to minimal generators
of $M/tM, M/sM,$ clearly hypotheses (i), (ii) are inherited by the
quotient $M/N$, so by induction we have
$$M/N\simeq (n-1)A/(t).$$ Now if $Q\not\subset (t)$, then $t$ is
regular (i.e. multiplication by it is injective) on $N=A/Q$ and
consequently the kernel of multiplication by $t$ on $M$ itself
maps isomorphically to $M/N=(n-1)A/(t).$ Therefore we get a
splitting
$$M\simeq (n-1)A/(t)\oplus N.$$
So $N$, as quotient of $M$, inherits properties (i), (ii) and
hence by the case $n=1$ already considered, we have $N\simeq
A/(t),$ contradiction. \par Therefore we may assume $Q\subset (t)$
hence, as $Q$ is primary, $Q=(t^m)$ for some $m\geq 1.$ Since $s$
is regular on $M/N$, we get an injection
$$N/sN\simeq A/(t^m, s)\hookrightarrow M/sM\simeq nA/(s,t)$$
which clearly forces $m=1.$ Now since $t$ kills $N$ and $M/N$ and
$N$ is saturated, being generated by a minimal generator $u$ of
$M$, it follows that $t$ kills $M$ so by (i), $M\simeq
nA/(t)$.\qed\enddemo

\proclaim{Corollary 8.3} If Proposition 8.1 holds for all $a\leq
2n-2$ then it holds for all $a$.\endproclaim\demo{proof}By
induction, suppose that $a>2n-2$ and that the Proposition holds
for all $a'<a.$ Consider a general curve of the form
$$C=C_{n-1}\cup_q C_{a-n+1}\cup_p C_b.$$ Consider a 2-parameter
smoothing of this curve parametrized by $s,t$ where $s=0$ (resp.
$t=0$) is the locus where $q$ (resp. $p$) remains singular.
Consider a general smoothing of the appropriate linear system, and
let $\N$ be the relative normalized normal bundle and set
$$N_0=N\otimes L\inv\otimes\O((-2)\cup (-r+1)\cup (-1))$$
(i.e. $\N$ twisted by a line bundle with the appropriate degrees
on the respective components of $C$). Now the results of \S6 show
that as $C_a$ specializes to $C_{a-n+1}\cup C_{n-1}$, the
cohomology of $\N$ remains constant: in fact $N_{C_{a-n+1}\cup
C_{n-1}}$ splits as a direct sum of line bundles in such a way
that the splitting deforms, summand by summand, to a splitting of
$N_{C_a}$ (cf. Remark 5.4.1). This implies that the cohomology of
$N_0$ is constant along the locus $t=0$ where $p$ remains a node
(in this locus the general curve is of the form $C_a\cup C_b$ and
the special one is $C_{n-1}\cup C_{a-n+1}\cup C_b$ ). By
induction, Proposition 8.1 holds for $C_{a-n+1}\cup C_b,$ and then
it follows easily that the restriction of the cohomology of $N_0$
(i.e. $R^1\pi_*(N_0)$, whose formation clearly commutes with
base-change, being an $R^{\text{top}}\pi_*$) on the locus $s=0$
where $q$ remains a node is annihilated by $t$, since the
corresponding assertion holds for the smoothing of $C_{a-n+1}\cup
C_b$. Applying the Lemma and restricting to a slice $s=\epsilon$
gives a smoothing as desired for Proposition 8.1. This proves the
Corollary.
\enddemo

We are therefore reduced to considering the case $a\leq 2n-2$
which implies $d\leq 3n-3$ and therefore $r\leq 5.$

Chopping (8.3) into short exacts, we see that Proposition 8.1 will
hold provided the cokernel of $R^1(\pi_*\phi)$ is of the form
$$\text{coker}(R^1(\phi)\simeq bk(0)\oplus
((n-1)(r+1)-(2d-2))\O_{B,0}\tag8.4$$ (note that if (8.4) holds
then $\pi_*(N_0)=0$ for the given smoothing $W$). Now
$R^1(\pi_*\phi)$ is a map of free $\O_B$-modules, and it is easy
to see that (8.4) holds iff the dual $\mu=R^1\pi_*(\phi)^*$
satisfies
$$\text{coker}(\mu)\simeq bk(0).\tag8.5$$
Now set $$V_{c,d}=\pi_*(\O(c,d)).$$Then written out explicitly,
using standard duality identifications and Lemma 7.1, $\mu$ comes
out as the map
$$W\otimes V_{r,0}\to V_{a+r-1,b}\oplus V_{a+r,
b-1},$$
$$\chi\otimes g\mapsto (g\nabla_1(\chi), g\nabla_2(\chi))$$
with $(\nabla_1,\nabla_2)$ as in Lemma 7.1. We must show a choice
of $W$ for which the image of $\mu$ is of the form $sz_1,...,sz_b,
z_{b+1},...,z_{2(d+r)}$ where $(z.)$ is a basis for the target of
$\mu$ and $s$ is the parameter on $B$.\par Now let us introduce
some convenient notation. Set
$$e=(a,b), e[m]=(a-m,b+m).$$
Then it is easy to see that
$$\nabla(u^i)=(\nabla_1(u^i), \nabla_2(u^i))=e[i]u^i,$$
$$\nabla(v^i)=e[-i]v^i.$$ Now we claim that the desired property
(8.5) will hold for $W$ with basis $\chi_0,...,\chi_n$ as follows:
$\chi_0=1, \chi_1=u$ and for each $2\leq i\leq n,$ either
$$\chi_i=u^{\nu_i}$$ (which we call a $u$-move), or
$$\chi_i=u^{\nu_i}+v^{\rho_i}$$ (which we call a $v$-move),
where the choice of
the exponents $\nu_i, \rho_i$ is to be specified. Note that a
basis for the source of $\mu$ is given by
$$\chi_i\otimes u^j, i=0,...,n, j=0,...,r,$$
and we view the $\mu(\chi_i\otimes u^j)$ as arranged naturally in
a $(r+1)\times (n+1)$ matrix $\Theta$ and as a totally ordered
set, column by column. We denote by $\Theta_i$ the submatrix of
$\Theta$ containing the first $i+1$ columns. One choice of basis
for the target of $\mu$ is given by
$$e[i(j)]u^j, e[i'(j)]u^j,j=0,...,a+r-1, e[a]u^{a+r},$$
$$e[k(j)]v^j, e[k'(j)]v^j,j=1,...,b-1, e[-b]v^{b},$$
as long as $i(j)\neq i'(j), k(j)\neq k'(j), \forall j$ (which
implies that $e[i(j)], e[i'(j)]$ (resp. $e[k(j)], e[k'(j)]$)) are
linearly independent. Our aim is to devise a 'winning strategy'
for the choice of $\nu_i$ and $\rho_i$'s, which by definition
means that we have
$$\text{im}(\mu)=\text{span}(e,k^2u,...,k^2u^{a+r-1},e[a]u^{a+r},
e[-1]v,...,e[-b]v^b,\tag8.6$$
$$se[-1],se[-2]v,...,se[-b]v^{b-1})=: \Xi $$
(span meaning over $\O_B$), which would clearly imply (8.5). That
$\im(\mu)\subseteq\Xi$ is trivial and will emerge in the ensuing
discussion.

To this end, set
$$U=V_{r,0}, W_i=\span(\chi_0,...,\chi_i).$$
Note that $\mu(W_1\otimes U)$ contains
$$\span(e[i],e[i-1])u^i=\span(k^2u^i), i=1,...,r$$
but only a 1-dimensional subspace of $k^2u^0, k^2u^{r+1}$.
Accordingly, we set
$$x_1=r, z_1=1$$
Generally, we will say that an entry $\theta_{ij}=u^k+*v^*$ of
$\Theta$ or a power $u^k$ is {\it{doubled}} if $$\mu(W_i\otimes
U)\supset k^2u^k;\tag8.7$$ to show (8.7) holds it is of course
sufficient-- and in practice necessary-- to show that
$e[j]u^k,e[j']u^k$ occur in $\Theta_i, j\neq j'.$ Now suppose that
$$\nu_0=0,\nu_1=2,...,\nu_{i-1}$$
have been defined, where $i\leq n$, in such a way that for each
$c<i,$, the set of $j$ such that $\mu(W_c\otimes U)$ contains
(resp. meets nontrivially) $k^2u^j$ forms an integer interval
$[1,x_c]$ (resp. $[0,y_c]$); indeed $y_c=\max(\nu_d:d\leq c)+r.$
We will specify a uniquely determined way of choosing $\nu_{i+1}$,
in case of a $u$-move, or $(\nu_{i+1},\rho_{i+1})$, in case of a
$v$-move. This will reduce our problem of proving (8.6) to that of
devising a suitable strategy of choosing at each stage a $u$ or
$v$ move, where the total number of $u$ (resp. $v$) moves is
$n-1-b$ (resp. $b$). We call such a sequence of moves {\it
allowable.}\par For a $u$-move, define
$$\nu'_{i}=x_{i-1}+1,$$
$$\nu_{i}=\min(\nu'_i,a),\tag8.8$$
$$ x_i=\min(y_{i-1}, \nu_i+r), y_i:=\max(\nu_d:d\leq i)+r.$$
We call the case $\nu'_i>a$ an {\it{overboard}} case. It is then
clear that $\mu(W_i\otimes U)$ contains (resp. meets nontrivially)
$k^2u^j$ iff $j\in [1,x_i]$ (resp. $j\in [0, y_i]),$ and that a
$u$-move essentially always enlarges the rank of $\mu(W\otimes U)$
by $r+1$: more precisely
$$\rk(\mu(W_i\otimes U))=\min(\rk(\mu(W_{i-1}\otimes U)+r+1,
2(d+r)).\tag8.9$$\par In the case of a $v$-move, define
$$\nu'_{i}=x_{i-1}-1,$$
$$\nu_i=\min(\nu'_i,a),\tag 8.10$$ and $x_i,y_i$ as in the $u$ case.
Here again an {\it{overboard}} case is where $\nu'_i>a.$ As for
$\rho_i$, we simply define
$$\rho_i=\max(\rho_c:c<i)+1,$$
where the $\max$ is over those $c$ so that $\rho_c$ is defined
(i.e. so that the $c$th move is a $v$-move), or 0 if  no such $c$
exist. Thus the $\rho$'s which are defined simply take the values
$1,...,b.$ Note that for a $v$-move, the $i$-th column of $\Theta$
takes the form (transposed)
$$(e_0u^{x_{i-1}-1}+e_1v^{\rho_i},e_0u^{x_{i-1}}+se_1v^{\rho_i-1},...)$$
where $e_0=e[x_{i-1}-1], e_1=e[-\rho_i].$ For the first two
entries in the vector, the $u$ term in in $\mu(W_{i-1}\otimes U)$
by definition, and consequently
$$e_1v^{\rho_i},se_1v^{\rho_i-1}\in \mu(W_{i}\otimes U).$$ From
the third entry on, the $v$ parts, which take the form
$$s^2e_1v^{\rho_i-2},...,s^{\rho_i}e_1,...,s^{\rho_i}e_1u^{r-\rho_i}$$
in case $\rho_i\leq r$, or
$$s^2e_1v^{\rho_i-2},...,s^{r}e_1v^{\rho_i-r}$$
in case $\rho_i\geq r$, are clearly in $\mu(W_{i-1}\otimes U)$, so
we conclude that
$$e_0u^{x_{i-1}+1},...,e_0u^{x_{i-1}+r-1}\in \mu(W_{i}\otimes U).$$
It is clear in any case that as long as $x_{i-1}+1\leq a$ ($u$
move) or $x_{i-1}-1\leq a$ ($v$ move), the entries of $\Theta_i$
form part of a basis of $\Xi.$ \par Now it is easy to see from the
definition that we always have
$$\nu_i\leq x_i\leq \nu_i+r.$$
Consequently,
$$\nu_i\geq \nu_{i-1},$$ except in the one case where the $i$th column is
a $v$-move and $x_{i-1}=\nu_{i-1}$ (we call this an
{\it{exceptional}} $v$-move), in which we have
$$\nu_i=\nu_{i-1}-1, x_i=\nu_i+r.$$ Accordingly, the top $u$ power
occurring in $\Theta$ through the $i$th column always occurs in
the $i$th column itself, except if the $i$th column is an
exceptional $v$-move, in which the top $u$ power occurs in the
$(i-1)$st column.\par It will be useful to introduce the notion of
$i$th {\it{level}} $z_i$, defined as follows: if
$\nu_i\geq\nu_{i-1},$ $z_i$ is the number of non-doubled elements
in the $i$th column, i.e. $z_i=\nu_i+r-x_i;$ if
$\nu_i=\nu_{i-1}-1,$ set $z_i=-1.$ Levels transform nicely under
$u$ and $v$ moves: define
$$\phi(x)=r+1-|x|, \psi(x)=r-1-|x|.$$
Then after a $u$ move in the $i$th column, $z_i=\phi(z_{i-1})$,
while after a $v$ move in the $i$th column, $z_i=\psi(z_{i-1})$.
Note that $\phi,\psi$ are essentially reflections, i.e.
$$\phi^2(x)=\psi^2(x)=|x|.$$ Now a key observation is the
following: \proclaim{Lemma 8.4} Suppose $\Theta$ is constructed by
an allowable sequence of moves and $z_n=1.$ Then (8.6)
holds\endproclaim To see this, suppose first there are no
overboards. Then the fact that $z_n=1$ implies that the powers
$u,...,u^{\nu_n+r-1}$ are all doubled and having $b$ many
$v$-moves ensures that
$$e[-i]v^i, se[-i-1]v^i \in\im(\mu), i=0,...,b-1 .$$ Also
$e[\nu_n]u^{\nu_n+r},e[-b]v^b\in\im(\mu),$ obviously. Since
$$(r+1)(n+1)\geq 2(d+r)=2a+2b+2r,$$ this is only possible if
$\nu_n=a$, so (8.6) holds. Now if an overboard occurs at step $i$,
it is easy to check that $\nu_j\geq a\ \forall j\geq i$ and using
$\nu_n=a$ we can conclude as above.\qed\par Recalling that $z_1=1$
automatically, we are thus reduced to finding an {\it{allowable
word}} $w$ in  $\phi,\psi$, i.e. one containing $b$ many $\psi$'s
and $n-1-b$ many $\phi$'s, so that $w(1)=1.$ This is easiest if
$n$ is odd: ideed in this case it suffices to set
$$w=\psi^b\phi^{n-1-b}, b\  \text{even},$$
$$w=\phi^{n-2-b}\psi^b\phi=\phi^{n-2-b}\psi^{b-1}(\psi\phi), b\ \text{odd}$$
(using $\psi\phi(1)=-1$). \par Henceforth we assume
that $n$ is
even. Recall from Corollary 8.3 that we may assume $2\leq r\leq
5$. Since we are assuming $(d,n)$ is perfect, we have
$$(r+1)(n-1)=2d-2,$$ therefore $r$ must be odd, i.e. $r=3$ or $5$.
Suppose first that
$$r=3.$$ This implies
$$\psi(\pm 1)=1.$$ So it suffices to take
$$w=\phi^{n-1-b}\psi^b, b\ \text{odd},$$
$$w=\phi^{n-2-b}\psi^b\phi=\phi^{n-2-b}\psi^{b-1}(\psi\phi), b\ \text{even}.$$
Next, suppose $$r=5.$$ Note that if $b=n-1$ then $a\geq n-1$; if
strict inequality holds, then $C_{a,b}$ is almost balanced by the
results of \S6, while if equality holds then $r=4$ (and $(d,n)$
cannot be perfect). Therefore we may assume $$2\leq b\leq n-2.$$
Then note
$$\psi\phi\psi(1)=1$$ so we can take
$$w=\phi^{n-2-b}\psi^{b-1}\phi\psi, b\ {\text {even}},$$
$$w=\phi^{n-b-3}\psi^{b-2}(\phi\psi)^2, b\ {\text {odd}}$$
(note that $b\leq n-2$, $b$ odd, $n$ even forces $n-b\geq 3.$)\par
\remark{Remark 8.4} Our only use in the above argument of the
assumption that $(d,n)$ is perfect  was to avoid the cases
$r=2,4$. However, with a fairly short additional argument that we
now sketch, the proof can be extended to cover those cases as
well. This involves taking advantage of the 'slack' due to the
fact that $(r+1)(n+1)>2(d+r)$, i.e. the source of $\mu$ has larger
rank than its asserted image. Indeed, define a map $\mu'$
analogous to $\mu$ but with the exponents $\nu'_i$ in place of
$\nu_i.$ Then, it is easy to check that we have a winning strategy
(for $\mu$) if we can arrange for $\mu'$ to have its  $n$th level
$z'_n$ be either 0, 1 or 2. The latter can be achieved with the
following choices.\par\noindent For $r=2$ (so $\phi(1)=2),
\psi(1)=0) $:
$$b\ \even : \ w=\phi^{n-1-b}\psi^b \Rightarrow z'_n=2$$ $$b\ \odd:
\ w=\psi^b\phi^{n-1-b}\Rightarrow z'_n=0.$$ For $r=4$ (so
$\phi(2)=3, \psi(1)=2$):
$$b\ \odd: \ w=\phi^{n-1-b}\psi^b\Rightarrow z'_n=2,$$ $$b\ \even: \
w=\psi\phi^{n-1-b}\psi^{b-1}\Rightarrow z'_n=1.$$

\endremark
\remark{Remark 8.5} Clearly condition (8.6), which is sufficient
for the conclusion of Proposition (8.1) to hold, depends on $W$
only modulo $s^2$, i.e. depends only on the angle $C_{a,b}$ and
its first-order deformation corresponding to $W\mod
s^2$.\endremark

 \comment

Note that then
$$x_{i+1}=\nu_i+r.\tag\?$$ Proceeding in this way, we define
$\nu_i, x_i, i\leq n-b$. Note that $\nu_{i+2}=\nu_i+r+1, i+2\leq
n-b$, so in fact,
$$\nu_{2k}=k(r+1), 2k\leq n-b,\tag\?$$
$$\nu_{2k+1}=1+k(r+1), 2k+1\leq n-b,\tag\?$$ while the $x_i$ can
be read off from \?, which yields
$$x_{2k}=\nu_{2k},
x_{2k+1}=\nu_{2k+1}+r-1.$$
 It is then easy to see inductively that $\mu|_{W_{n-b}\otimes
U}$ is injective fibrewise (i.e. $\mod s$), that its image is
generated by 'monomials' of the form $zu^j, z\in k^2$, and that
$\mu(W_{i}\otimes U)$ contains $k^2u^j$ for $j\in [1,x_i]$ but
only a 1-dimensional subspace of $k^2u^j$ generated by $e[k]u^j,
k\leq\mu_i $ for $j=0$ or $j\in [x_i+1,\nu_i+r].$\par Next, define
$$\nu_{n-b+1}=x_{n-b}-1.$$
Then by definition $\mu(W_{n-b+1}\otimes U)$ contains
$$\chi_{n-b+1}=e[n-b+1]u^{n-b+1}+e[-1]v,$$ and since the first summand
is already in $\mu(W_{n-b}\otimes U)$, we conclude
$$e[-1]v\in \mu(W_{n-b+1}\otimes U).$$ Next,
$u\chi_{n-b+1}\in \mu(W_{n-b+1}\otimes U) $ yields similarly
$$se[-1]\in\mu(W_{n-b+1}\otimes U).$$
Continuing similarly, we see that
$$u^{x_{n-b}+1}e[n-b+1]+se[-1]u\in\mu(W_{n-b+1}\otimes
U),$$ hence $u^{x_{n-b}+1}e[n-b+1]\in\mu(W_{n-b+1}\otimes U)$.
Since $\mu(W_{n-b}\otimes U)$ already contains $u^{x_{n-b}+1}e[j]$
for some $j<n-b+1$, we see that $\mu(W_{n-b+1}\otimes U)$ contains
$k^2u^{x_{n-b}+1}$ so there is a number $x_{n-b+1}>x_{n-b}$ such
that $\mu(W_{n-b+1}\otimes U)$ contains $k^2u^j$ iff $j\in
[1,x_{n-b+1}]$
Indeed it is easy to see that in
case $n-b$ is even, so $x_{n-b}=\nu_{n-b}=(n-b)(r+1)/2,$ we have
$$\nu_{n-b+1}=\nu_{n-b}-1, x_{n-b+1}=\nu_{n-b}+r-1,$$
and  $\mu(W_{n-b+1}\otimes U)$ contains $e[n-b]u^{\nu_{n-b}+r}$,
while if $n-b$ is odd we have
$$\nu_{n-b+1}=\nu_{n-b}+r-2, x_{n-b+1}=\nu_{n-b}+r$$
and $\mu(W_{n-b+1}\otimes U)$ contains $e[n-b+1]u^{\nu_{n-b}+j},
j\in [r+1, 2r-2]$. For general $n-b< i<n$, if $\nu_i,x_i$ have
been defined, we set
$$\nu_{i+1}=\min(x_i-1, a)$$ and let $x_{i+1}$ be defined by the
property that $\mu(W_{i+1}\otimes U)$ contains $k^2u^j$ iff $j\in
[1,x_{i+1}]$.
In fact, the pattern of $\nu_i,x_i$ becomes quite regular with
period 2 after $i=n-b+1$ or $i=n-b+2$ and one can easily verify
the following. \par {Case $n-b$ even:}\par set
$\nu=\nu_{n-b}=(n-b)(r+1)/2.$ Then
$$\nu_{n-b+2}=\nu+r-2, x_{n-b+2}=\nu+r,$$
$$e[\nu+r-2]u^j\in\mu(W_{n-b+2}\otimes U), j\in [\nu+r+1, \nu+2r=2]
$$ and for $i\geq n-b+2,$
$$\nu_{i+1}=\min(\nu_i+1,a),x_{i+1}=\nu_{i+1}+r-1,
e[\nu_{i+1}]u^{\nu_{i+1}+r}
\in\mu(W_{i+1}\otimes U), \  i\ \text{even};$$
$$\nu_{i+1}=\min(\nu_i+r-2,a), x_{i+1}=\nu_{i+1}+2,$$
$$e[\nu_{i+1}]u^j\in\mu(W_{i+1}\otimes U), j\in
[\nu_{i+1}+3,\nu_{i+1}+r],
 \ i\ \text{odd}.$$
Case $n-b$ odd:\par For all $i\geq n-b,$
$$\nu_{i+1}=\min(\nu_i+r-2,a), x_{i+1}=\nu_{i+1}+2,$$$$
e[\nu_{i+1}]u^j\in\mu(W_{i+1}\otimes U), j\in [\nu_{i+1}+3,
\nu_{i+1}+r], \ i\ \text{odd},$$
$$\nu_{i+1}=\min(\nu_i+1,a), x_{i+1}=\nu_{i+1}+r-1,$$$$
e[\nu_{i+1}]u^{\nu_{i+1}+r}\in\mu(W_{i+1}\otimes U), \ i\
\text{even}.$$ In any case, we can check as before that for all
$i\geq n-b+2,$$$e[-i+n-b]v^{i-n+b},
se[-i+n-b]v^{i+n-b-1}\in\mu(W_i\otimes U).$$ Now, by definition of
$r$, we have $$(r+1)(n+1)\geq 2d+2r=2a+2r+2b.$$ From this, is is
easy to see that at the end, $\mu(W\otimes U)$ contains $k^2v^j$
for all $j\in [1,a]$, as well as $e[-j]v^j, se[-j-1]v^j$ for all
$j\in [0,b-1]$ and $e[-b]v^b,$ which proves our assertion \?
\endcomment
\enddemo
\comment

 ******************

 To begin with, the restriction of $N_0$ on
$X_0$ is given by \? :
$$N_0|_{X_0}=(d-n)\O_{X_0}(0\cup 0)\oplus (a_1-1)\O_{X_0}(0\cup(-2))$$
$$\oplus (b_1-1)\O_{X_0}((-2)\cup
0)\oplus\O_{X_0}((-1)\cup(-1)):=N_{0,0}\oplus N_{0,1}\oplus
N_{0,2}\oplus N_{0,3}. \tag\?$$

We let
$$N_0\to N_{0,2}\oplus N_{0,3}$$ denote the corresponding
surjection and $N_1$ its kernel, which is thus an elementary
reduction of $N_0.$ Also let $\N_1=N_1(2\tau_1+\tau_2)$ be the
corresponding elementary reduction of $\N$. Since
$\pi_*(N_{0,2}\oplus N_{0,3})=0$, while $R^1\pi_*((N_{0,2}\oplus
N_{0,3})$ is a $k(0)$-vector space of dimension $n-a,$ the obvious
cohomology sequence yields
$$\l_0(R^1\pi_*(N_1))=\l_0(R^1\pi_*(N_0))-(n-a).\tag\?$$
At this juncture, a critical point is the determination of the
restriction of $N_1$ on $X_0$:

\proclaim{Lemma \?} We have
$$ N_1|_{X_0}\simeq (d-n-1)\O_{X_0}(0\cup 0)\oplus N_{0,1}\oplus
N_{0,2}\oplus \O_{X_0}((-1)\cup 0)\oplus\O_{X_0}(0\cup
(-1))\tag\?$$
\endproclaim
\demo{proof} By way of explanation and interpretation, it is easy
to see a priori that \? is equivalent to the conjunction of the
following assertions:\newline (\?.1) $N_1|_{X_i}\simeq N_0|_{X_i},
i=1,2$;\newline (\?.2) The HN filtrations of $N_1|_{X_1}$ and
$N_1|_{X_2}$ are in {\it{general position}} (not just relative
general position) at the node $p$.\newline /*****************

 Note
that by the construction of elementary reductions, the subsheaf
$HN_1(N_0|_{X_i})\simeq (a-1)\O_{X_i}$ survives in $N_1$, $i=1,2$
hence
$$HN_1(N_1|_{X_1})\simeq HN_1(N_0|_{X_1})\simeq (a-1)\O_{X_1}$$
and similarly for $X_2$, and since the $HN_1$'s are already in
general position for $N_0$, clearly they remain so for $N_1.$

Given (\?.1) and \?, the only issue  (2) is equivalent to\newline
(2') $HN_1(N_1|_{X_1})(p)\cap HN_1(N_1|_{X_2})(p)\subsetneq
HN_1(N_0|_{X_1})(p)\cap HN_1(N_0|_{X_2})(p).$ ****************/
\par
Now consider the exact sequence on $X$
$$0\to\omega\inv\to T\otimes L\inv\overset\psi\to \to \N\to q\to 0\tag\?$$
where $\omega=\omega_{X/B}$ and $q$ is a skyscraper sheaf at the
node $p=([1,0],[1,0])$. Let's work in affine coordinates
$u=U_1/U_0, v=V_1/V_0$ near $p$. Let $(\chi.)$ be an adapted basis
for $W$ ordered so that in affine coordinates, $\chi_0=1, \chi_i$
for $i=1,...,a_1$ has order $i$ on $X_2,$ while $\chi_{n_1-i}$ for
$i=1,...,b_1$ has order $i$ on $X_1.$ Then as $T\otimes L\inv$ is
the pullback of tautological quotient bundle on $\P^n$, it is the
quotient of $W$ by the rank-1 subsheaf generated by
$$f=(\chi_0,...,\chi_n)=$$
$$(1, u*,...,u^{a_1}*, u^{a_1+1}*+v^b*,..., u^a*+v^{b_1}*,
v^{b_1}*,...,v*)$$ where the $*$ denote units. Since the relative
dualizing sheaf is generated by $(du\wedge dv)ds\inv$, the
relative tangent sheaf $\omega\inv$ is generated by $ds(du\wedge
dv)\inv$, i.e. by the vector field $\partial=u\partial/\partial
u+v\partial/\partial v,$ therefore the image $\N'$ of $\psi$ is
the quotient of $\pi^*W$ by the subsheaf generated by $f$ and
$$f'=\partial f=(0,\partial\chi_1,...,\partial\chi_n)$$
i.e. the row-space $D$ of the Wronski matrix
$$W[f]=\left [\matrix f\\f'\endmatrix\right ].$$
Note that
$$u^i|\partial\chi_i,\forall 1\leq i\leq a_1,\tag\?$$
$$v^i|\partial\chi_{n+1-1},\forall 1\leq i\leq b_1.\tag\?$$  Thus we have an
exact sequence
$$ 0\to D\to \pi^*W\to \N'\to 0.\tag\?$$ Note that $\N|_{X_1}$
admits a natural filtration (which coincides with its HN
filtration, and which obviously splits for degree reasons) with
quotients
$$\N_{C_a/\P^a}=(a-1)\O_{X_1}(2), \O_{X_1}(1), (n-a-1)\O_{X_1}$$
and similarly for $X_2$. The image of $\N'|_{X_1}$ in $\N|_{X_1}$
i.e. $\N'|_{X_1}$ mod torsion, is just
$N_{C_a/\P^n}=(a-1)\O(2)\oplus (n-a)\O$, where the second factor
is just the restriction of $\N_{\P^a/\P^n}$, hence one of the $\O$
factors, viz. the one corresponding to $T_pC_b \mod \P^a$ maps to
the $\O(1)$ factor in $\N_{C_a/\P^a}$ (which is well-defined $\mod
\P^a$), the map vanishing at $p$. For a suitable local basis of
$\N$, the map $\pi^*W\to\N$ takes the standard basis vector $e_1$
to $(v,0...0)$ and $e_{n-b_1+1}$ to $(u,0...0)$. Also note that
the basic reduction sequence \? yields in the case of $\N_1$
$$ 0\to \O(1\cup 1)\oplus (b_1-1)\O(0\cup 2)\to \N_1|_{X_0}\to
(d-n)\O(2\cup 2)\oplus (a_1-1)\O(2\cup 0)\to 0\tag\?  $$

Restricting this further on $X_1,$ we get
$$0\to \O_{X_1}(1)\oplus (b_1-1)\O_{X_1}\to \N_1|_{X_1}\to
(a-1)\O_{X_1}(2)\to 0\tag\?$$
Viewing \? as an extension, it splits naturally in two components,
involving $\O_{X_1}(1)$ and $(b_1-1)\O_{X_1}$ respectively. We are
claiming that this extension splits. This is obvious as far as the
$\O_{X_1}(1)$ component, so the point is about the other
component, involving $\O_{X_1}(2)$ and $\O_{X_1}$.

\par
Now let $W_1\subset W$ be the subsheaf generated by
$e_0,...,e_{a}, se_{a+1},...,se_n$, so that $\pi^*W_1$ is the
elementary reduction of $\pi^*W$ corresponding to the directions
normal to $\P^a$, and is thus compatible with $\N_1$ in the sense
that there is an induced map
$$\pi^*W_1\to \N_1\tag\?.$$
Since the quotient $N_0\to q$ factors through the quotient $N_0\to
N_{0,2}\oplus N_{0,3}$ defining $N_1,$ clearly the induced map \?
is surjective. There is a rational homomorphism of $D$ to
$\pi^*W_1$, defined off $X_2,$ given by the matrix obtained by
dividing the last $b_1$ columns of $W[f]$ by $s$, i.e.
$$W_1[f]=\left [\matrix
1&\chi_1&...&\chi_a&\chi_{a+1}/s&...&\chi_{n-1}/s&\chi_n/s\\
0&\partial\chi_1&...&\partial\chi_{a}&\partial\chi_{a+1}/s&...
&\partial\chi_{n-1}/s&\partial\chi_n/s\endmatrix\right ]=$$
/*********
$$\left [\matrix
1&u+...&...&u^a+v^{b_1}&(v^{b_1-1}+...)/u&...&(v/u+...)&(1/u+...)\\
0&u+...&...&au^{a}+b_1v^{b_1}+...&((b_1-1)v^{b_1-1}+...)/u&...&(v/u+...)
&(1/u+...)\endmatrix\right ].\tag\?$$ ************/

The poles indicate that $D\not\subset\pi^*W_1$, but if we set
$D_1=D(X_1-X_2)$ although $D$ certainly contains a full-rank
subsheaf contained in $\pi^*W_1$, viz. $D_1:=D\cap \pi^*W_1.$
Restricting on $X_1$ and putting things together, we get an exact
diagram
$$\matrix 0\to&D_1|_{X_1}&\to&\pi^*W_1&\to&\N_1&\to&q&\to 0\\
              &\downarrow&&\downarrow&&\downarrow&&\downarrow\\
          0\to&D|_{X_1}&\to&\pi^*W&\to&\N&\to&q&\to 0\\
          &&&\downarrow&&\downarrow&&\\
          &&&b_1\O_{X_1}&=&b_1\O_{X_1}&&\endmatrix\tag\?$$
By \?, last $(b_1-1)$ columns of $W_1[f]$ are divisible by $v$
(albeit with denominator $u$), hence the map $D_1|_{X_1}\to
(b_1-1)\O_{X_1}$ given by these  columns  vanishes generically,
hence vanishes, so we get a map $\N_1|_{X_1}\to (b_1-1)\O_{X_1}$
which splits the $\O$ component of the extension \?. Thus the
extension splits, as claimed, so that
$$\N_1|_{X_1}\simeq \N_{X_1}\tag\?$$
and likewise for $X_2$ (of course, the isomorphism in question is
{\it not} induced by the natural inclusion $\N_1\to\N$). This, of
course, implies that
$$N_0|_{X_i}\simeq N_1|_{X_i}, i=1,2$$
which is the first assertion of the lemma.\par It is useful to
note at this point that this assertion is true, with essentially
the same proof, for an arbitrary smoothing of $C_{a,b}$: the point
is that the normalized normal bundle restricted on $C_a$ consists
of an 'internal' part, i.e. $(a-1)\O(2)$ and an external part
$\O(1)\oplus (n-a-1)\O$ and the $(n-a-1)\O$ part corresponds to
elements of the linear system vanishing on $C_a$ or $X_1$ to order
$>1.$\par

To continue, let's return to \? and view it as an extension. The
extension class (viewed as an element of $H^1$) naturally break up
into 4 components,
$$o_1\in (d-n)H^1(\O((-1)\cup(-1)), o_2\in
(a_1-1)H^1(\O((-1)\cup 1))=0,$$
$$o_3\in (b_1-1)(d-n)H^1(\O((-2)\cup 0)), o_4\in
(b_1-1)(a_1-1)\O((-2)\cup 2))\tag\?$$ We have just proven that the
extension splits on $X_1$, which implies that the image of $o_3,
o_4$ in $(b_1-1)H^1(\O_{X_1}(-2))$ (resp.
$(b_1-1)(a_1-1)H^1(\O_{X_1}(-2))$) is zero. This clearly implies
that
$$o_3=0, o_4=0.$$
In light of all this, it is now a simple formal exercise to see
that the general position assertion \? now is equivalent to
$$o_1\neq 0\tag\?$$ (thus the extension \? splits on each
component $X_i$ but not on $X_1\cup X_2$). The proof of \?
proceeds in 2 steps: first, use specialization and projection to
reduce to the 'smallest' case, i.e. $a=b=2, n=3.$; then for this
case, use an explicit calculation.\par Assume $a\leq b.$ To begin
with, we specialize the curve $C_{a,b}$ to a curve where $C_a$ and
$C_b$ meet in $a-2$ points $p_1,...,p_{a-2}$ besides $p$, but are
otherwise general. Clearly this can be done, now just to $C_{a,b}$
but to a given 'smoothing' of it (which will be a smoothing at $p$
but preserve the crossings at the other nodes. The new curve will
be nondegenerate if $b+3\geq n$, otherwise it will span a
$\P^{b+3}.$ Then pick $b-a$ further general points
$p_{a-1},...,p_{b-2}$ on $C_b$ and project from $p_1,...,p_{b-2}$
to $\P^{n-b+2}$. We obtain a (general) smoothing of a conic
$C_{2,2}\subset\P^{n-b+2}$ We shall assume $n=b+1$, i.e.
$n-b+2=3$, otherwise $n-b+2\geq 4$ and the argument is simpler.
Now the projection induces a map on normalized normal bundles and
twists thereof, which extends to a map on their elementary
reductions. It follows from this that to prove \? for a general
smoothing of our original $C_{a,b}$ it suffices to prove it for a
general smoothing of $C_{2,2}\subset\P^3$ and since the assertion
involves an open property, it suffices to prove it for one
smoothing.\par We shall work in affine coordinates $(x,y,z)$ on
$\A^3,$ leaving it to the reader  to check that things extend
appropriately to $\P^3.$ Consider the family given parametrically
by
$$(u+u^2, u+v, -v+v^2)$$
that is a smoothing of the union of the two conics
$$C_2: x=y+y^2, z=0; C'_2: z=-y+y^2, x=0.$$
The planes of the conics meet in the $y-$axis and this meets $C_2$
in the origin $p$ and $p_1=(0,-1,0)$ (coming from $u=-1, v=0$),
and meets $C'_2$ in $p$ and $p_2=(0,1,0)$ (coming from $u=0,
v=1.$) Set
$$h=xz-s(1-y^2+x+z+s), q_1=y+y^2-x+z, q_2= z+y-y^2.$$
All these are polynomials vanishing on the family. Moreover it is
easy to see that $h$ yields the unique section of
$N^*(2L-p_1-p_2)$ on $X_0$ - note we know at this point that
$$N^*(2L-p_1-p_2)|_{X_0}=\O(0\cup 0)\oplus\O((-1)\cup(-1)).$$
Moreover, as sections of $N^*(2L)$ (though obviously not as
polynomials), $q_1$ and $q_2$ glue together to a section $q$
which, together with $h$, realizes the splitting
$$N^*(2L)|_{X_0}=\O(1\cup 1)\oplus\O.$$
Now our assertion \? is equivalent to the statement that after
performing the elementary reduction on $N^*(2L)$ corresponding to
the quotient $\O_{X_0}$, thus obtaining $N_1,$ say, the section
$h$, which survives as a section of $N_1|_{X_0},$ does {\it not}
survive as section of $N_1|_{X_0}(-p_1-p_2)$, i.e. there is no
constant $a$ such that $h+asq$ 'vanishes', in the appropriate
sense, on $p_1$ and $p_2$. To check this, note that local
generators for $N^*$ at $p_i$ are given by $q_i$ and
$$g_1=z+s(1-v)/u, g_2=x-s(1+u)/v$$
respectively for $i=1,2$. If our assertion is false, there is a
constant $a$ such that
$$h+asq_i\in\m_{p_i, X_0}(g_i, q_i), i=1,2\tag\?$$
For $i=1,$ we expand this in terms of the local parameter $(y+1)$
centered at $p_1$, and read off the leading term, which yields
$a=1$; for $i=2,$ doing likewise with $(y-1)$ at $p_2$ yields
$a=-2$. This contradiction finally proves the assertion \?.

\enddemo

\endcomment
\ss \subheading{9. Enumerative results}\ss We are now able to
state and prove our main enumerative results on rational curves
whose normal bundles are not almost balanced. This will be done
mainly by combining the cohomological computations of \S8 with the
(Grothendieck) Riemann-Roch formula, using as well some known
qualitative and enumerative results concerning a generic
'incidence pencil' of rational curves, i.e. a 1-parameter family
defined by incidence to a generic collection of rational curves,
which are summarized in the Appendix, whose notation and results
we shall be using freely.\par Thus let $$\pi:X\to B$$ be a generic
incidence pencil as in the Appendix, i.e. a smooth model of a
generic 1-parameter family of rational curves of degree $d$ in
$\P^n$ incident to a generic collection $(A.)$ of linear spaces.
Let
$$f:X\to \P^n$$ be the natural map, and set $L=f^*\O(1).$ We choose
an ordering on the set of components $X_1,X_2$ of each reducible
fibre of $\pi$, so that $$a:=L.X_1\geq b:=L.X_2,$$ and we call
such a fibre of {\it{bidegree}} $(a,b)$. We assume henceforth that
$(d,n)$ is a perfect pair and set, as in \S8
$$r=\frac{2d-2}{n-1}-1.$$
 A {\it{twisting divisor}} $D$ on $X$ is
by definition an integral divisor on $X$ such that\par (i) for a
fibre $F$ of $\pi,$ $D.F=r+2;$\par (ii) for a reducible fibre
$X_1\cup X_2$ of bidegree $(a,b)$, $$D.X_1=r+1, D.X_2=1.$$ It is
clear that twisting divisors exist: a specific choice is given by
$$D=(r+1)s_1-\sum\limits_{X_2\in\F_1}X_2-
r\sum\limits_{X_1\in\F_1}X_1,\tag9.1$$ where $\F_1$ is the set of
fibre components not meeting $s_1.$\par We fix a twisting divisor
$D$ and let $\N$ be the relative normalized normal bundle of
$X/B$, and set
$$G=\N(-D).$$ Then the restriction of $G$ on a smooth fibre $X_b$ of
$\pi$ has degree $-(n-1)$. In fact for the generic $b$, $f(X_b)$
is a generic rational curve of degree $d$ in $\P^n$ and hence by
Theorem 6.1, we have $$G|_{X_b}\simeq (n-1)\O(-1).\tag9.2$$ The
smooth fibres $X_b$ for which (9.2) does not hold are precisely
those for which
$$h^1(G|_{X_b})>0.$$ Accordingly, we call the corresponding curve
$C=f(X_b)$ a {\it{normally jumping rational curve}} of
multiplicity $h=h^1(G|_{X_b}),$ a number which for the generic
pencil depends only on $C$ itself. \proclaim{Theorem 9.1} With the
above notations, the number of normally jumping curves in the
pencil $B(a.)$, weighted according to multiplicity, is given by
$$J^\perp_{n,d}(a.)=\hskip 2in\tag9.3$$$$L(L+K_{X/B})
-D(\frac{n-3}{2}K_{X/B}+\frac{n-1}{2}L) -2N_d^{\text{red}}(a.)
-\sum (b-1)N_{a,b}(a.),$$ the summation being over all $a+b=d,
b\leq a, 1<b<n.$
\endproclaim\demo{proof} Applying (-1) times Grothendieck-Riemann-Roch
\cite{F} to the sheaf $G$ and the mapping $\pi$, we get an
equation, one side of which is the total length of $R^1\pi_*(G)$,
since by the above discussion $\pi_*(G)=0$ and $R^1\pi_*(G)$ is of
finite length. By Proposition 8.1, the length in question is the
sum of the sought-for weighted number of normally jumping curves
and the $\sum (b-1)N_{a,b}(a.)$ term from (9.3). To complete the
proof it suffices to evaluate the other side of -GRR, which is
routine. Briefly, the general formula of \cite{F} yields
$$-\int\limits_X
(((n-1)1_X+c_1(G)+\frac{1}{2}(c_1^2-2c_2)(G))(1_X-\frac{1}{2}K_X
+\chi(\O_X)[{\text{pt}}]))_2\tag9.4$$ where [pt] is a point and
$(\cdot)_2$ denotes the part in degree 2. Note that
$K_X=\omega+(2g-2)F$ where $\omega=K_{X/B}$, $F$ is a fibre of
$\pi$ and $g$ is the genus of $B$. Now the Chern classes of $\N$
can be computed from the exact sequence (8.1), yielding
$$c_1(\N)=2L+\omega,$$
$$c_2(\N)=3L(L+\omega)+\omega^2-\sigma$$
where $\sigma=N_d^{\text{red}}(a.)$ is the number of singular
points of (fibres of) $\pi$, which by (A15) equals $-\omega^2.$
Therefore
$$(c_1^2-2c_2)(\N)=-L^2-2L\omega-3\sigma.$$ Then standard Chern
calculus yields
$$\frac{-1}{2}(c_1^2-2c_2)(G)=L^2+L\omega-\frac{3}{2}\sigma+
(2L+\omega)D)-\frac{n-1}{2}D^2.\tag9.5$$ For the other two product
terms of (9.4) , note that because $X$ is a blown-up
$\P^1$-bundle, we have $\chi(\O_X)=\chi(\O_B)=1-g$ and by (9.2),
$c_1(G).F=-(n-1),$ hence
$$-\chi(\O_X)+\frac{1}{2}K_Xc_1(G)=\frac{1}{2}\omega(2L+\omega-(n-1)D).
\tag9.6$$ Then summing (9.5) and (9.6) and subtracting $\sum
(b-1)N_{a,b}(a.)$ yields (9.3).\qed
\enddemo\ls
 \heading Appendix: rational curves in $\P^n$\footnote{
This Appendix is reproduced with a few modifications from
\cite{R1} for the reader's convenience}\endheading

The purpose of this appendix is to review some notations and
results, both qualitative and enumerative, about rational curves
in $\P^n$ that are used in the statement and proof of the main
result. Proofs and further details may be found in \cite{R2, R3,
R4} and references therein.

We begin by reviewing some qualitative results about families of
rational curves in $\P^n$, especially for $n\geq 3.$ See [R2][R3]
[R4] and references therein for details and proofs. In what
follows we denote by $\bar{V}_d$ the closure in the Chow variety
of the locus of irreducible nonsingular rational curves of degree
$d$ in $\P^n, n\geq 3$, with the scheme structure as closure, i.e.
the reduced structure (recall that the Chow form of a reduced
1-cycle $Z$ is just the hypersurface in $G(1,\P^3)$ consisting of
all linear spaces meeting $Z$). Thus $\V_d$ is irreducible reduced
of dimension
$$\dim (\V_d)=(n+1)d+n-3.$$
Let $$A_1, \ldots, A_k\subset\P^n$$ be a generic collection of
linear subspaces of respective codimensions $a_1, \ldots, a_k,$ $
1\leq a_i \leq n$.
 We denote
by
$$
B= B_d = B_d (a_{\cdot}) = B_d (A_{\cdot})
$$
the normalization of the locus (with reduced structure)
$$
\{(C, P_1, \ldots, P_k ) \ : \ C \in \bar{V}_d , P_i  \in C \cap
A_i, i =1 , \ldots, k \},
$$ and refer to it as a (generic) {\it{incidence family}} or
{\it{incidence subvariety}} of $\bar{V}_d.$ If all $a_i>1$ then
this locus
 is also the normalization of its projection to $\bar{V}_d$,
i.e. the locus of degree-$d$ rational curves (and their
specializations) meeting $A_1, \ldots, A_k$.
 We have
$$
\dim B = (n+1)d+n-3  - \sum (a_i - 1) .\tag A1
$$
When $\dim B=0$ we set
$$N_d(a.)=\deg (B). \tag A2$$
Note that $N_d(1,a.)=dN_d(a.)$, which allows us to reduce the
computation of the general $N_d(a.)$ to the case where all
$a_i>1,$ in which case we will say the condition-vector $(a.)$ is
{\it{proper}}. The number $k$ of $a_i$ such that $a_i>1$ is called
the {\it length} of the condition-vector $(a.)$. Whenever
$b=\dim(B)\geq 0,$ it is convenient to set
$$N'_d(a.)=N_d(a.,b+1)\tag A3$$
and note that this is the degree in $\P^n$ of the locus swept out
by the curves in $B(a.)$.

The numbers $N_d$ and $N_d(a.)$, first computed in general by
Kontsevich and Manin (see for instance [FP] and references
therein),
 are computed in [R2],[R3] by an elementary method, reviewed below,
based on recursion on $d$ and $k$.

Now suppose $B=B(a.)$ is such that $\dim B=1$ and $(a.)$ is proper
and let
$$
\pi : X \to B
$$
be the normalization of the tautological family of rational
curves, and $$f: X \to \P^n$$ the natural map. We call $B$ or
$X/B$ a (generic) {\it{incidence percil}}.  The following
summarizes results from [R2][R3][R4] : \proclaim{Theorem A0}(i) X
is smooth .\par (ii) Each fibre $C$ of $\pi$ is either\par (a) a
$\P^1$ on which $f$ is either an immersion with at most one
exception which maps to a cusp ($n=2$) or an embedding ($n>2$); or
\par (b)  a pair of
$\P^1$'s meeting transversely once, on which $f$ is an immersion
with nodal image ($n=2$) or an embedding ($n>2$); or
\par (c)  if $n=3$, a $\P^1$ on which $f$ is a degree-1
immersion such that $f(\P^1)$ has a unique singular point which is
an ordinary node.\par (iii) If $n>2$ then $\bar{V}_{d,n}$ is
smooth along the image $\bar{B}$ of $B$, and $\bar{B}$ is smooth
except, in case some $a_i=2,$ for ordinary nodes corresponding to
curves meeting some $A_i$ of codimension 2 twice. If $n=2$ then
$\bar{V}_{d,n}$ is smooth in codimension 1 except for a cusp along
the cuspidal locus and normal crossings along the reducible locus,
and $\bar{B}$ has the singularities induced from $\bar{V}_{d,n}$
plus ordinary nodes corresponding to curves with a node at some
$A_i$, and no other singularities.
\endproclaim

Next, we review some of the enumerative apparatus introduced in
[R3][R4] to study $X/B$. Set
$$m_i=m_i(a.)=-s_i^2, i=1,...,k.\tag A4 $$
Note that if $a_i=a_j$ then $m_i=m_j.$
 It is shown in
[R2] [R3][R4] that these numbers can all be computed recursively
in terms of data of lower degree $d$ and lower length $k$. For
$n\geq 2,$ note that
$$
s_i . s_j = N_d(...,a_i+a_j,...,\hat a_j,...), i \neq j.
$$

 Also, letting
$R_\l$ denote the sum of all fibre components not meeting $s_\l$ ,
we have
$$
s_j \cdot R_\l = \sum N'_{d_1} (a_i:i\in I ) N'_{d_2}
(a_i:i\not\in I ). \tag A6
$$
the summations being over all $d_1 + d_2 = d$ and all index-sets
$I$ with $j\in I, \l\not\in I$, so all these numbers may be
considered known. Then we have
$$m_i=\frac{1}{2}(s_i.R_j+s_i.R_p-s_j.R_p)-s_i.s_j
-s_i.s_p+s_j.s_p $$ for any distinct $i,j,p,$ and the RHS here is
an expression of lower degree and/or length, hence may be
considered known.
\par
Next, set $$L=f^*(\O(1)),$$ and note that
$$L^2=N_d(2,a.),\  L.s_i=N_d(a_1,...,a_{i}+1,...)
,\ i=1,...,k$$ (in particular, $L.s_i=0$ if $a_i=n.$) We computed
in [R3] that, for any $i$,
$$
L \sim d s_i - \sum\limits_{F \in \F_{i}} \deg (F) F + ( N_d (a_1
,\ldots, a_i+1, \ldots) + dm_i(a.)) F_0
$$
where $F_0$ is the class of a complete fibre and $\F_i$ is the set
of fibre components not meeting $s_i$. Consequently we have
$$
N_d(2,a_1, \ldots) = 2d N_d (a_1 + 1, a_2, \ldots) + d^2 m_1(a.) -
\sum\limits_{F\in \F_{1}(a.)} (\deg F)^2 \tag A7
$$
and clearly the RHS is a lower degree/length expression, so all
the $N_d(2,\ldots)$ are known. We also have for $n>2$ that
$$
N_d(a_1,a_2+1,...) - N_d(a_1+1,a_2,...)=$$$$
dN_d(a_1+a_2,...)-\sum\limits_{F \in (\F_{1}-\F_{2})(a.)} (\deg F)
+ N_d(a_1+1,a_2,...)+dm_1(a.) \tag A8
$$
and again the RHS here is 'known', hence so is the LHS, which
allows us to 'shift weight' between the $a_i$'s till one of them
becomes equal to 2, so we may apply (A7), and thus compute all of
the $N_d(a.)$'s.\par Next, it is easy to see as in [R3] that

$$
L.R_j = \sum d_2N'_{d_1} (a_i:i\in I ) N_{d_2} (
 a_i:i\not\in I ),\ n\geq 2
\tag A10 $$ the summation for $n>2$ being over all $d_1 + d_2 =
d,$ and all index-sets $I$ such that $j\in I.$\par Finally, the
relative canonical class $K_{X/B}=K_X-\pi^*(K_B)$ was computed in
[R3] as
$$K_{X/B}=-2s_i-m_iF+R_i\tag A11$$
for any $i$. Note that $-R_i^2$ equals the number $\sigma$ of
reducible (equivalently, singular) fibres in the family $X/B$, a
number we denote by $N_d^{\text{red}}(a_.)$, and which is easily
computable by recursion, namely let
$$N_{d_1,d_2}=
\binom{3d-2}{3d_1-1}d_1d_2 N_{d_1}N_{d_2}, n=2,$$
$$N_{d_1,d_2}(a_.)=\sum\limits_I N'_{d_1}(a_i:i\in I)N'_{d_2}
(a_i:i\not\in I), n\geq 2,\tag A12$$ where the latter sum extends
over all index-sets $I$. Then
$$N_d^{\text{red}}(a_.)=\sum\limits_{d_1+d_2=d}N_{d_1,d_2}(a_.).\tag
A13$$

From this we compute easily that
$$L.K_{X/B}= -2N_d(...a_i+1...)-dm_i+L.R_i,\tag A14$$
$$K_{X/B}^2=-\Ndr(a_.).\tag A15$$


\Refs\widestnumber\key{EV2}

\ref\key EV1\by D. Eisenbud, A. Van de Ven\paper On the normal
bundle of smooth space curves\jour Math. Ann. \vol 256 \yr
1981\pages 453-463\endref
 \ref\key EV2\bysame\paper On the variety of smooth
rational space curves with given degree and normal bundle\jour
Invent. math. \vol 67 \yr 1982\pages 89-100\endref \ref\key F\by
W. Fulton\book Intersection theory \publ Springer\yr 1984\endref
 \ref\key GS\by
F. Ghione, G. Sacchiero\paper Normal bundles of rational curves
 in $\P\sp 3.$
\jour Manuscr. Math. \vol 33\pages 111-128\yr 1980\endref \ref\key
H\by R. Hartshorne\book Algebraic Geometry\publ Springer \yr
1977\endref\ref\key OSS\by  Ch. Okonek, M. Schneider, H. Spindler
\book Vector bundles on projective spaces\publ Birkh\"auser \yr
1980\endref\ref\key Ram\by L. Ramella\paper La stratification du
sch\'ema de Hilbert des courbes rationelles de $\P^n$ par le
fibr\'e tangent restreint\jour Comptes Rendus Acad. Sci. Paris
S\'er. I\vol 311 \yr 1990\pages 181-184\endref \ref\key{R1} \by Z.
Ran \paper Geometry on nodal curves (preprint)\endref \ref\key
R2\bysame\paper Bend, break and count \jour Isr. J. Math\vol 111
\yr 1999 \pages 109-124\endref \ref\key R3\bysame \paper Bend,
break and count II \jour Math. Proc. Camb. Phil . Soc. \vol 127\yr
1999\pages 7-12
\endref
\ref\key R4\bysame\paper On the variety of rational space
curve\jour Isr. J. Math \vol 122\yr 2001\pages 359-370\endref
\ref\key R5\bysame\paper The degree of the divisor of jumping
rational curves\jour Quart. J. Math.\yr 2001
 \pages 1-18\endref
 \ref\key R6\bysame\paper Enumerative geometry of divisorial
families of rational curves \finalinfo eprint math.AG/0205090,
updated version at $\underline{www.math.ucr.edu/\tilde{\ }
ziv/papers/ geonodal.pdf}$\endref
\endRefs
\enddocument